\documentclass{amsart}

\usepackage{amssymb}
\usepackage{amsmath}
\usepackage{amscd}
\usepackage{ascmac}

\usepackage{graphicx,amsfonts}
\usepackage{wrapfig}

\usepackage{bm}

\theoremstyle{plain}
\newtheorem{thm}{Theorem}[section]
\newtheorem{cor}[thm]{Corollary}
\newtheorem{lem}[thm]{Lemma}
\newtheorem{prop}[thm]{Proposition}

\theoremstyle{definition}
\newtheorem{defn}[thm]{Definition}

\theoremstyle{remark}
\newtheorem{rem}[thm]{Remark}

\begin{document}

\title{On sections of hyperelliptic Lefschetz fibrations}
\author[S.~Tanaka]{Shunsuke Tanaka}
\address{Department of Mathematics, Graduate School of Science,
Osaka University, Toyonaka, Osaka 560-0043, Japan}
\email{qqx33r5d@almond.ocn.ne.jp}
\date{January 20, 2012; MSC 2000: primary 57N13, secondary 20F34}

\maketitle

\begin{abstract} 
We construct a relation among right-handed Dehn twists in the mapping class group 
of a compact oriented surface of genus $g$ with $4g+4$ boundary components. 
This relation gives an explicit topological description of $4g+4$ disjoint $(-1)$-sections of 
a hyperelliptic Lefschetz fibration of genus $g$ on the manifold 
$\mathbb{CP}^2\#(4g+5)\overline{\mathbb{CP}}^2$. 
\end{abstract}

\section{Introduction}

Lefschetz fibrations relate the topology of symplectic $4$-manifolds 
to the combinatorics on relations in Dehn twist generators of mapping class groups 
of surfaces. 
It is well-known that a Lefschetz fibration of genus $1$ on the manifold 
$E(1)=\mathbb{CP}^2\# 9\overline{\mathbb{CP}}^2$ constructed by blowing up nine 
intersections of two generic cubics in $\mathbb{CP}^2$ has twelve singular fibers 
and nine disjoint $(-1)$-sections. 
Korkmaz and Ozbagci \cite{KO} constructed a relation among right-handed 
Dehn twists in the mapping class group of a torus with nine boundary components 
to locate a set of nine disjoint $(-1)$-sections in a Kirby diagram of $E(1)$. 
It is also known to algebraic geometers that a hyperelliptic Lefschetz fibration of 
genus $g$ on the manifold $X_g=\mathbb{CP}^2\# (4g+5)\overline{\mathbb{CP}}^2$ 
has $8g+4$ singular fibers and $4g+4$ disjoint $(-1)$-sections for $g\geq 2$ 
(cf. \cite[Sect. 3]{SS}, see also \cite[Remark 1.1]{KK}). 

In this paper we construct a relation among right-handed Dehn twists 
in the mapping class group of a compact oriented surface of genus $g$ 
with $4g+4$ boundary components to locate a set of 
$4g+4$ disjoint $(-1)$-sections in a Kirby diagram of $X_g$. 
In the case $g=2$, our relation is considered as 
an improvement of Onaran's relations \cite{onaran} 
in mapping class groups of surfaces of genus two 
with at most eight boundary components.

In Section 2 we review basic relations in mapping class groups 
and produce two relations on a torus with eight boundary components. 
Combining these relations, 
we construct a new relation on a surface of genus $g$ with $4g+4$ boundary components 
in Section 3. 
In Section 4 we apply the relation to visualize $4g+4$ disjoint $(-1)$-sections in a 
Kirby diagram of a hyperelliptic Lefschetz fibration of genus $g$. 

Several results of this paper are based on master's thesis of the author 
at Osaka University in February, 2011.

\section{Building blocks}

In this section we review basic relations in mapping class groups 
and produce two relations on a torus with boundary 
both used in the next section.

\subsection{Basic relations in mapping class groups}

Let $\Sigma_{g,r}$ be a compact oriented surface of genus $g$ 
with $r$ boundary components 
and ${\rm Diff}_+\Sigma_{g,r}$ the group of 
orientation-preserving diffeomorphisms of $\Sigma_{g,r}$ fixing the boundary 
$\partial \Sigma_{g,r}$ pointwise equipped with $C^{\infty}$-topology. 
The group $\pi_0({\rm Diff}_+\Sigma_{g,r})$ of path-components of 
${\rm Diff}_+\Sigma_{g,r}$ 
is called the {\it mapping class group} of $\Sigma_{g,r}$ and we denote it by 
$\mathcal{M}_{g,r}$. 
We denote by $\mathcal{F}_{g,r}$ the free group 
generated by all isotopy classes $\mathcal{S}_{g,r}$ of 
simple closed curves in the interior of  $\Sigma_{g,r}$. There is a natural 
epimorphism $\varpi:\mathcal{F}_{g,r}\rightarrow \mathcal{M}_{g,r}$ which sends 
(the isotopy class of) a simple closed curve $a$ in the interior of $\Sigma_{g,r}$ 
to the right-handed Dehn twist $t_a$ along $a$. 
We set $\mathcal{R}_{g,r}:={\rm Ker}\; \varpi$. 

A word in the generators $\mathcal{S}_{g,r}$ is called {\it positive} if 
it includes no negative exponents. 
We put ${}_W(c):=t_{a_r}^{\varepsilon_r}\cdots 
t_{a_1}^{\varepsilon_1}(c)\in \mathcal{S}_{g,r}$ for $c\in \mathcal{S}_{g,r}$ and 
$W=a_r^{\varepsilon_r}\cdots a_1^{\varepsilon_1}\in \mathcal{F}_{g,r} \; 
(a_1, \ldots ,a_r \in \mathcal{S}_{g,r}, \varepsilon_1,\ldots ,\varepsilon_r\in \{ \pm 1\})$. 
We often denote $a^{-1}$ by $\bar{a}$ for an element $a$ of $\mathcal{S}_{g,r}$. 
For two words $W_1,W_2\in\mathcal{F}_{g,r}$, we denote $W_1\equiv W_2$ if 
$\varpi(W_1)=\varpi(W_2)$. 

We recall definitions of basic relations in mapping class groups. 

\begin{figure}[h]
\includegraphics[width=3cm]{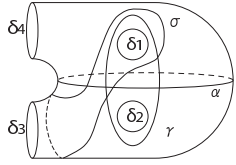}
\caption{\label{lantern} Lantern relation}
\end{figure}

\begin{figure}[h]
\includegraphics[width=8cm]{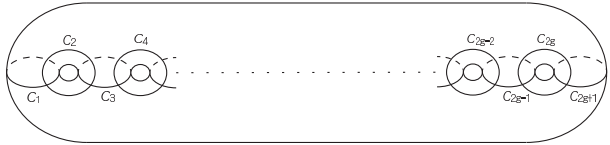}
\caption{\label{hyperelliptic} Hyperelliptic relation}
\end{figure}

\begin{defn}[cf. \cite{FM}]\label{basic}
(1) For disjoint simple closed curves $a,b$ in the interior of $\Sigma_{g,r}$, we have 
a relation $ab\equiv ba$ in $\mathcal{F}_{g,r}$ 
called a {\it commutativity relation}. 
A regular neighborhood of $a\cup b$ is the disjoint union of two annuli. 

(2) For simple closed curves $a,b$ in the interior of $\Sigma_{g,r}$ 
which intersect transversely at one point, 
we have a relation 
$aba\equiv bab$ in $\mathcal{F}_{g,r}$ 
called a {\it braid relation}. 
A regular neighborhood of $a\cup b$ is a torus with one boundary component. 

(3) For simple closed curves $\alpha, \sigma, \gamma, \delta_1,\delta_2,\delta_3,\delta_4$ 
in the interior of $\Sigma_{g,r}$ shown in Figure \ref{lantern}, 
we have a relation $\delta_1\delta_2\delta_3\delta_4\equiv \gamma\sigma\alpha$ 
in $\mathcal{F}_{g,r}$ called a {\it lantern relation}. 
The union of $\delta_1,\delta_2,\delta_3,\delta_4$ bounds a sphere with four boundary 
components in $\Sigma_{g,r}$. 

(4) An ordered $n$-tuple $(c_1,\ldots ,c_n)$ of simple closed 
curves in the interior of $\Sigma_{g,r}$ is called a {\it chain} of length $n$ 
if $c_i$ and $c_{i+1}$ intersect transversely at one point 
$(i=1,\ldots ,n-1)$ and other $c_i$ and $c_j$ never intersect. 
For a chain $(c_1,\ldots ,c_{2g+1})$ of length $2g+1$ on $\Sigma_{g,0}$, 
we have a relation $(c_1\cdots c_{2g+1}c_{2g+1}\cdots c_1)^2\equiv 1$ in 
$\mathcal{F}_{g,0}$ called a {\it hyperelliptic relation} (cf. Figure \ref{hyperelliptic}). 
\end{defn}

\vspace{0cm}

\begin{rem}\label{braid}
Let $a$ and $b$ be simple closed curves in the interior of $\Sigma_{g,r}$ and 
$c$ the simple closed curve $t_b(a)={}_b(a)$. Then we have the relation 
$c\equiv ba\bar{b}$ in $\mathcal{F}_{g,r}$. 
If $a$ and $b$ intersect transversely at one point, 
we have another relation $b\equiv ac\bar{a}$. 
This relation together with the relation 
$c\equiv ba\bar{b}$ yields a braid 
relation $aba\equiv bab$. 
\end{rem}

\subsection{Two relations on a torus with boundary}

In this subsection we construct two relations on a torus with eight boundary components. 
The first relation is the following. 

\begin{figure}[h]
\includegraphics[width=5cm]{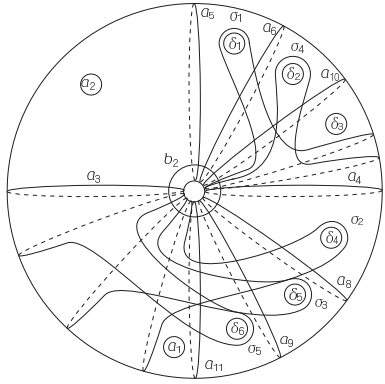}
\caption{\label{relA} Relation A}
\end{figure}

\begin{prop}[Relation A]\label{relationA}
For simple closed curves in the interior of $\Sigma_{1,8}$ shown in Figure \ref{relA}, 
we have the relation 
\[
a_1 a_{2} \delta_1 \delta_2 \delta_3 \delta_4 \delta_5 \delta_6
\equiv a_5a_4 b_2 \bar{a}_4 \sigma_1 \sigma_4 a_{10} 
a_3 b_2 \bar{a}_3 \sigma_2 a_5 a_3 a_8 b_2 \bar{a}_8 \bar{a}_3 \sigma_3 \sigma_5 a_{11}. 
\]
\end{prop}

We make use of the five-holed torus 
relation found by Korkmaz and Ozbagci \cite{KO} 
in order to prove Proposition \ref{relationA}.

\begin{figure}[h]
\includegraphics[width=10cm]{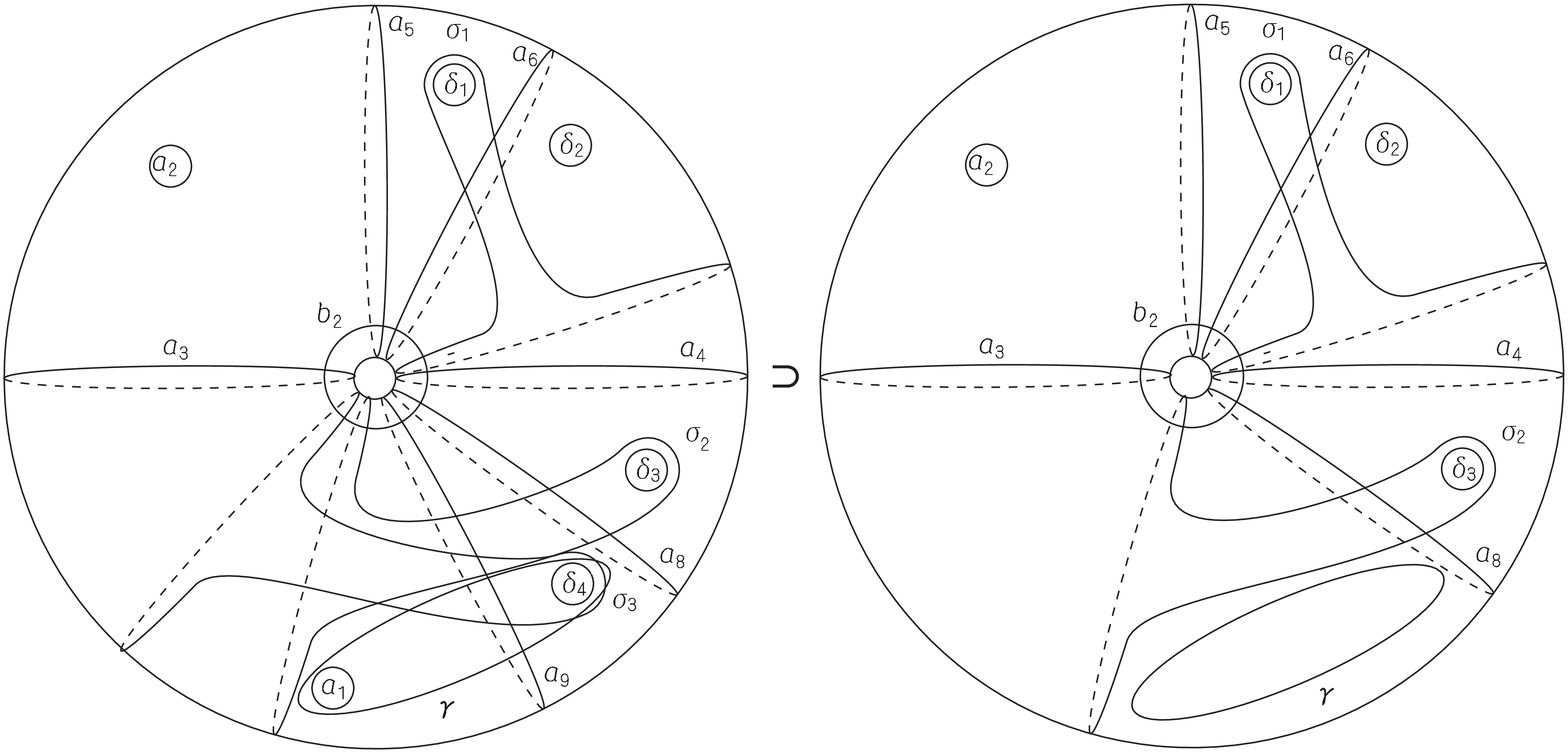}
\caption{\label{five} Five-holed torus relation}
\end{figure}

\begin{lem}[Korkmaz-Ozbagci \cite{KO}]\label{five-holed}
For simple closed curves in the interior of $\Sigma_{1,5}$ 
shown on the right in Figure \ref{five}, 
we have the relation 
\[
\delta_2\delta_1a_2\gamma\delta_3
\equiv a_5 b_2 a_3 a_4 a_5 b_2 \sigma_1 a_6 a_3 b_2 \sigma_2 a_8. 
\]
\end{lem}

%\vspace{-1.5cm}

\noindent
{\it Proof of Proposition \ref{relationA}}. 
Applying commutativity relations and conjugations to the five-holed torus relation 
in Lemma \ref{five-holed}, we obtain 
\begin{align*}
a_2\delta_1\delta_2\delta_3\gamma 
& \equiv a_5b_2a_3a_4a_5b_2\sigma_1a_6a_3b_2\sigma_2a_8 
\equiv a_8b_2\bar{a}_8a_8a_3a_4a_5b_2\sigma_1a_6a_3b_2\sigma_2a_5 \\ 
& \equiv a_8a_3a_4a_5b_2\sigma_1a_6a_3b_2\sigma_2a_5a_8b_2\bar{a}_8. 
\end{align*}
Multiplying both sides of this relation by $\bar{\gamma}$, we have 
\[
a_2\delta_1\delta_2\delta_3 
\equiv a_8a_3a_4a_5b_2\sigma_1a_6a_3b_2\sigma_2a_5a_8 b_2\bar{a}_8\bar{\gamma}. 
\]
We embed $\Sigma_{1,5}$ into $\Sigma_{1,6}$ and take simple closed curves 
$a_1,a_9,\delta_4,\sigma_3$ in the interior of $\Sigma_{1,6}$ 
shown in Figure \ref{five}. Then we have a lantern relation 
\[
\delta_4a_1a_3a_8 \equiv \gamma \sigma_3 a_{9}. 
\]
Combining these relations and applying commutativity relations, we obtain 
\begin{align*}
a_8a_3a_1a_2\delta_1\delta_2\delta_3\delta_4 
& \equiv a_8a_3a_4a_5b_2\sigma_1a_6a_3b_2\sigma_2a_5a_8 
b_2 \bar{a}_8\bar{\gamma}\gamma\sigma_3a_{9}     \\
& \equiv a_8a_3a_4a_5b_2\sigma_1a_6a_3b_2\sigma_2a_5a_8b_2\bar{a}_8\sigma_3a_{9}.     
\end{align*}
Multiplying both sides of this relation by $\bar{a}_3\bar{a}_8$, we have a relation 
\begin{equation}\tag{A1}
a_1a_2\delta_1\delta_2\delta_3\delta_4 
\equiv a_4a_5b_2\sigma_1a_6a_3b_2\sigma_2a_5a_8b_2\bar{a}_8\sigma_3 a_{9}       
\end{equation}
on $\Sigma_{1,6}$. 

\begin{figure}[h]
\includegraphics[width=10cm]{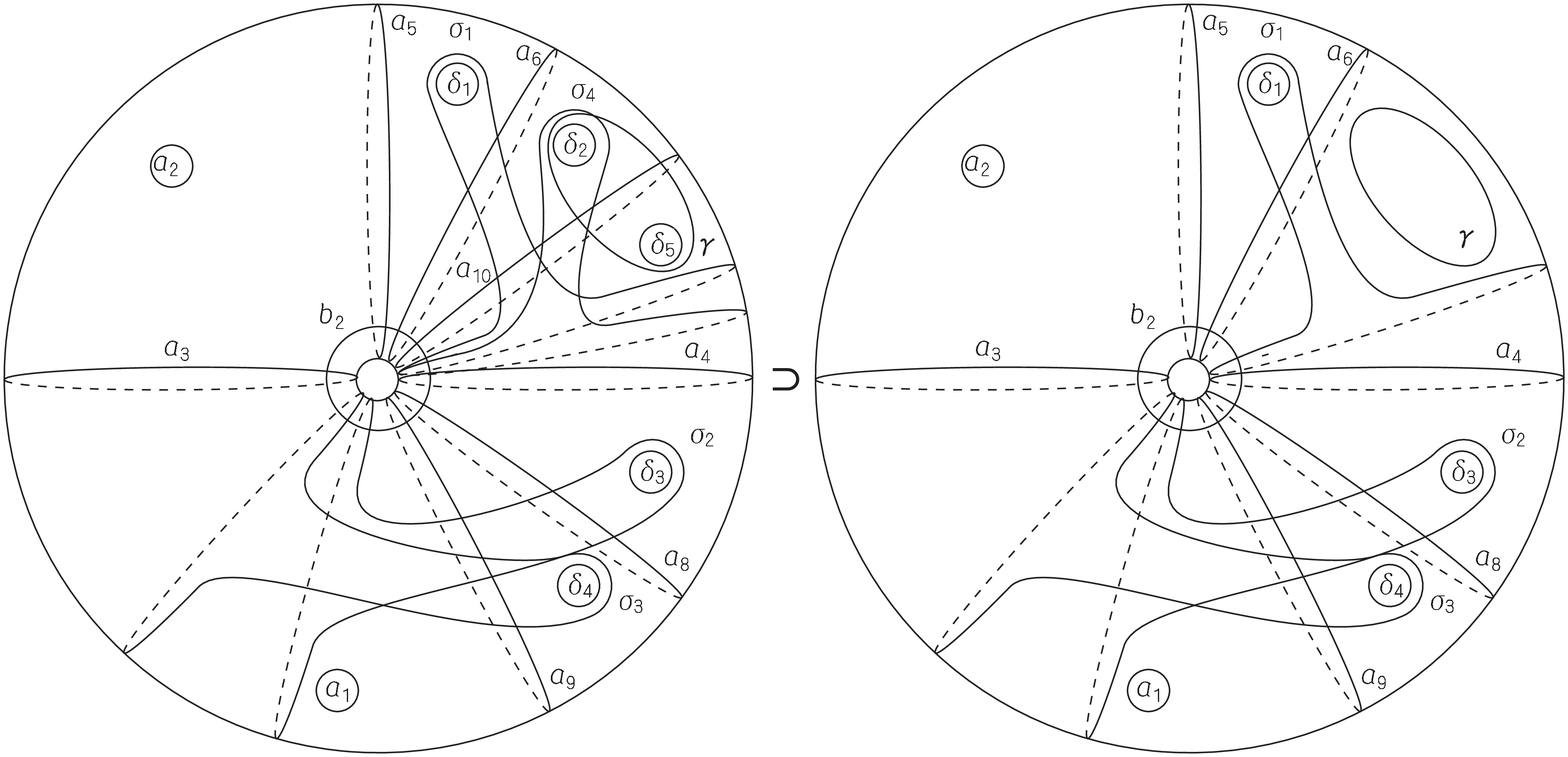}
\caption{\label{6to7} Embedding of $\Sigma_{1,6}$ into $\Sigma_{1,7}$ (I)}
\end{figure}

We change the name $\delta_2$ of a curve in the relation (A1) 
into $\gamma$ (shown on the right in Figure \ref{6to7}) 
and apply commutativity relations and conjugations to it to obtain 
\begin{align*}
a_1a_2\delta_1\delta_3\delta_4\gamma 
& \equiv a_4a_5b_2\sigma_1a_6a_3b_2\sigma_2a_5a_8 b_2\bar{a}_8\sigma_3a_9 \\
& \equiv a_5a_4 b_2\bar{a}_4a_4\sigma_1a_6a_3b_2\sigma_2a_5a_8b_2
\bar{a}_8\sigma_3a_{9}  \\ 
& \equiv a_4a_6a_3b_2\sigma_2a_5a_8b_2\bar{a}_8\sigma_3a_9a_5a_4b_2\bar{a}_4\sigma_1.    
\end{align*}
Multiplying both sides of this relation by $\bar{\gamma}$, we have 
\[
a_1a_2\delta_1\delta_3\delta_4 
\equiv a_4a_6a_3b_2\sigma_2a_5a_8b_2\bar{a}_8\sigma_3a_9a_5a_4b_2
\bar{a}_4\sigma_1 \bar{\gamma}. 
\]

\noindent
We embed $\Sigma_{1,6}$ into $\Sigma_{1,7}$ and take simple closed curves 
$a_{10},\delta_2,\delta_5,\sigma_4$ in the interior of $\Sigma_{1,7}$ 
shown in Figure \ref{6to7}. Then we have a lantern relation 
\[
\delta_2\delta_5a_4a_6 \equiv \gamma\sigma_4a_{10}. 
\]
Combining these relations and applying commutativity relations, we obtain 
\begin{align*}
a_4a_6a_1a_2\delta_1\delta_2\delta_3\delta_4\delta_5
& \equiv a_4a_6a_3b_2\sigma_2a_5a_8b_2\bar{a}_8\sigma_3a_9a_5a_4b_2\bar{a}_4\sigma_1
\bar{\gamma}\gamma\sigma_4a_{10} \\ 
& \equiv a_4a_6a_3b_2\sigma_2a_5a_8b_2\bar{a}_8\sigma_3a_{9}a_5a_4b_2\bar{a}_4
\sigma_1 \sigma_4 a_{10}. 
\end{align*}
Multiplying both sides of this relation by $\bar{a}_6\bar{a}_4$, we have a relation 
\begin{equation}\tag{A2}
a_1a_2\delta_1\delta_2\delta_3\delta_4\delta_5 
\equiv a_3b_2\sigma_2a_5a_8b_2\bar{a}_8\sigma_3a_9a_5a_4b_2\bar{a}_4
\sigma_1\sigma_4a_{10} 
\end{equation}
on $\Sigma_{1,7}$. 

We change the name $a_1$ of a curve in the relation (A2) 
into $\gamma$ (shown on the right in Figure 6) 
and apply commutativity relations and conjugations to it to obtain 
\begin{align*}
a_2\delta_1\delta_2\delta_3\delta_4\delta_5\gamma
& \equiv a_3b_2\sigma_2a_5a_8b_2\bar{a}_8\sigma_3a_9a_5a_4 b_2\bar{a}_4
\sigma_1\sigma_4a_{10} \\
& \equiv a_3b_2\bar{a}_3\sigma_2a_5a_3a_8b_2\bar{a}_8\sigma_3 
a_9a_5a_4 b_2\bar{a}_4\sigma_1\sigma_4a_{10} \\
& \equiv a_3b_2\bar{a}_3\sigma_2a_5a_3a_8b_2\bar{a}_8\bar{a}_3\sigma_3 
a_3a_9a_5a_4b_2\bar{a}_4\sigma_1\sigma_4a_{10} \\
& \equiv a_3a_9a_5a_4b_2\bar{a}_4\sigma_1\sigma_4a_{10} 
a_3b_2\bar{a}_3\sigma_2a_5a_3a_8b_2\bar{a}_8\bar{a}_3\sigma_3. 
\end{align*}
Multiplying both sides of this relation by $\bar{\gamma}$, we have 
\[
a_2\delta_1\delta_2\delta_3\delta_4\delta_5 
\equiv a_3a_9a_5a_4b_2\bar{a}_4\sigma_1\sigma_4a_{10} 
a_3b_2\bar{a}_3\sigma_2a_5a_3a_8b_2\bar{a}_8\bar{a}_3\sigma_3\bar{\gamma}. 
\]

\begin{figure}[h]
\includegraphics[width=10cm]{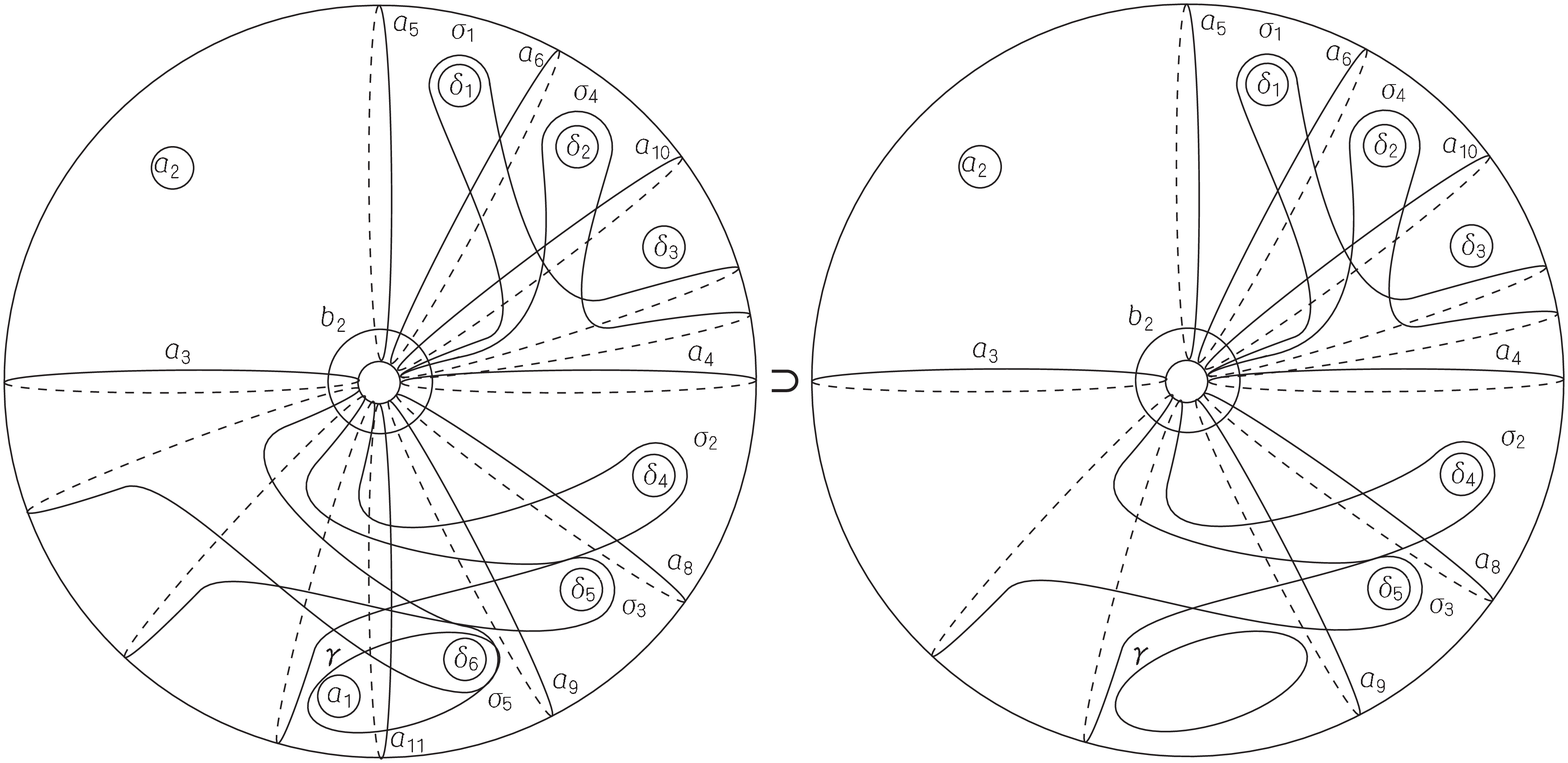}
\caption{\label{7to8} Embedding of $\Sigma_{1,7}$ into $\Sigma_{1,8}$ (I)}
\end{figure}

\noindent
We embed $\Sigma_{1,7}$ into $\Sigma_{1,8}$ and take simple closed curves 
$a_1,a_{11},\delta_6,\sigma_5$ in the interior of $\Sigma_{1,8}$ 
shown in Figure \ref{7to8}. Then we have a lantern relation 
\[
\delta_6a_1a_3a_9 \equiv \gamma\sigma_5a_{11}. 
\]
Combining these relations and applying commutativity relations, we obtain 
\begin{align*}
a_3a_9a_1a_2\delta_1\delta_2\delta_3\delta_4\delta_5\delta_6
& \equiv a_3a_9a_5a_4b_2\bar{a}_4\sigma_1\sigma_4a_{10} 
a_3b_2\bar{a}_3\sigma_2a_5a_3a_8b_2\bar{a}_8\bar{a}_3
\sigma_3\bar{\gamma}\gamma\sigma_5a_{11} \\
& \equiv a_3a_9a_5a_4b_2\bar{a}_4\sigma_1\sigma_4a_{10} 
a_3b_2\bar{a}_3\sigma_2a_5a_3a_8b_2\bar{a}_8\bar{a}_3\sigma_3\sigma_5a_{11}. 
\end{align*}
Multiplying both sides of this relation by $\bar{a}_9\bar{a}_3$, 
we finally obtain Relation A. This completes the proof of Proposition \ref{relationA}. 
$\square$

\medskip

The second relation constructed in this subsection is the following. 

\begin{figure}[h]
\includegraphics[width=5cm]{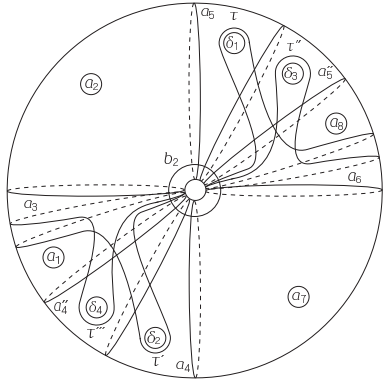}
\caption{\label{relB} Relation B}
\end{figure}

\begin{prop}[Relation B]\label{relationB}
For simple closed curves in the interior of $\Sigma_{1,8}$ shown in Figure \ref{relB}, 
we have the relation 
\[
a_1 a_2 a_7 a_8 \delta_{1} \delta_{2} \delta_{3} \delta_{4} 
\equiv a_4 a''_5 \bar{a}_6 b_2 a_6 a_3 b_2 \bar{a}_3 \tau' \tau''' 
a_5 a''_4 \bar{a}_3 b_2 a_3 a_6 b_2 \bar{a}_6 \tau \tau''. 
\]
\end{prop}

We make use of the four-holed torus 
relation found by Korkmaz and Ozbagci \cite{KO} 
in order to prove Proposition \ref{relationB}.

\begin{figure}[h]
\includegraphics[width=10cm]{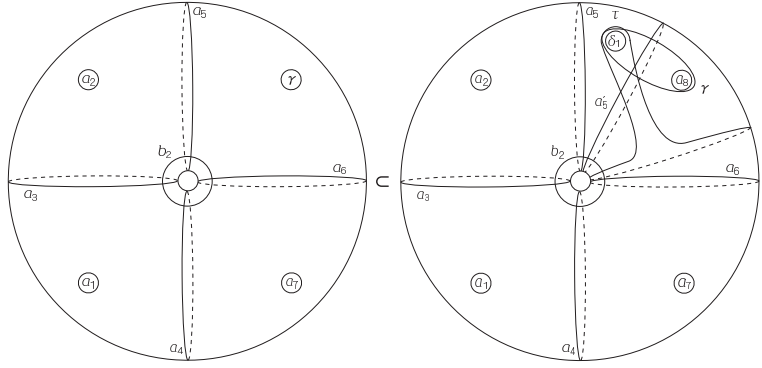}
\caption{\label{four} Four-holed torus relation}
\end{figure}

\begin{lem}[Korkmaz-Ozbagci \cite{KO}]\label{four-holed}
For simple closed curves in the interior of $\Sigma_{1,4}$ 
shown on the left in Figure \ref{four}, 
we have the relation 
\[
a_2a_1a_7\gamma \equiv (a_3a_6b_2a_4a_5b_2)^2. 
\]
\end{lem}

\vspace{0cm}

\noindent
{\it Proof of Proposition \ref{relationB}}. 
We consider the four-holed torus relation reviewed in Lemma \ref{four-holed}. 
We then embed $\Sigma_{1,4}$ into $\Sigma_{1,5}$ and take simple closed curves 
$a'_5,a_8,\delta_1,\tau$ in the interior of $\Sigma_{1,5}$ 
shown in Figure \ref{four}. Then we have a lantern relation 
\[
\delta_1a_8a_6a_5 \equiv  \gamma\tau a'_5. 
\]
Combining this relations with the four-holed torus relation, 
and applying commutativity relations and conjugations, we obtain a relation 
\begin{align*}\tag{B1}          
a_1a_2a_7a_8\delta_1 
& \equiv \bar{a}_5a_4a_5b_2a_3a_6b_2a_4a_5b_2a_3a_6b_2  
\bar{\gamma}\cdot\bar{a}_6\gamma\tau a'_5 \\    
& \equiv a_4b_2a_3a_6b_2a_4a_5b_2a_6a_3b_2\bar{a}_6\tau a'_5 \\   
& \equiv a_4a_5b_2a_6a_3b_2\bar{a}_6\tau a'_5a_4b_2a_3a_6b_2               
\end{align*}
on $\Sigma_{1,5}$. 

\begin{figure}[h]
\includegraphics[width=10cm]{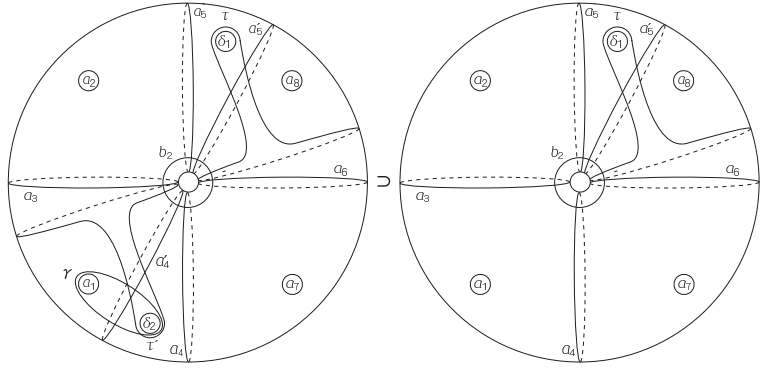}
\caption{\label{5to6} Embedding of $\Sigma_{1,5}$ into $\Sigma_{1,6}$}
\end{figure}

We change the name $a_1$ of a curve in the relation (B1) 
into $\gamma$ (shown on the right in Figure \ref{5to6}) to obtain 
\[
\gamma a_2a_7a_8\delta_1
\equiv a_4a_5b_2a_6a_3b_2\bar{a}_6\tau a'_5a_4b_2a_3a_6b_2. 
\]
We embed $\Sigma_{1,5}$ into $\Sigma_{1,6}$ and take simple closed curves 
$a_1,a'_4,\delta_2,\tau'$ in the interior of $\Sigma_{1,6}$ 
shown in Figure \ref{5to6}. Then we have a lantern relation 
\[
a_4a_3a_1\delta_2 \equiv  \gamma\tau' a'_4. 
\]
Combining these relations 
and applying commutativity relations and conjugations, we obtain a relation 
\begin{align*}\tag{B2}
a_1a_2a_7a_8\delta_1\delta_2 
& \equiv \bar{a}_4a_4a_5b_2a_6a_3b_2\bar{a}_6\tau a'_5a_4b_2a_3a_6b_2
\bar{\gamma} \cdot \bar{a}_3 \gamma \tau' a'_4 \\    
& \equiv a_5b_2a_6a_3b_2\bar{a}_6\tau a'_5a_4b_2a_6a_3b_2\bar{a}_3\tau' a'_4 \\    
& \equiv b_2a_6a_3b_2\bar{a}_3\tau' a'_4a_5b_2a_6a_3b_2\bar{a}_6\tau a'_5a_4.          
\end{align*}
on $\Sigma_{1,6}$. 

\begin{figure}[h]
\includegraphics[width=10cm]{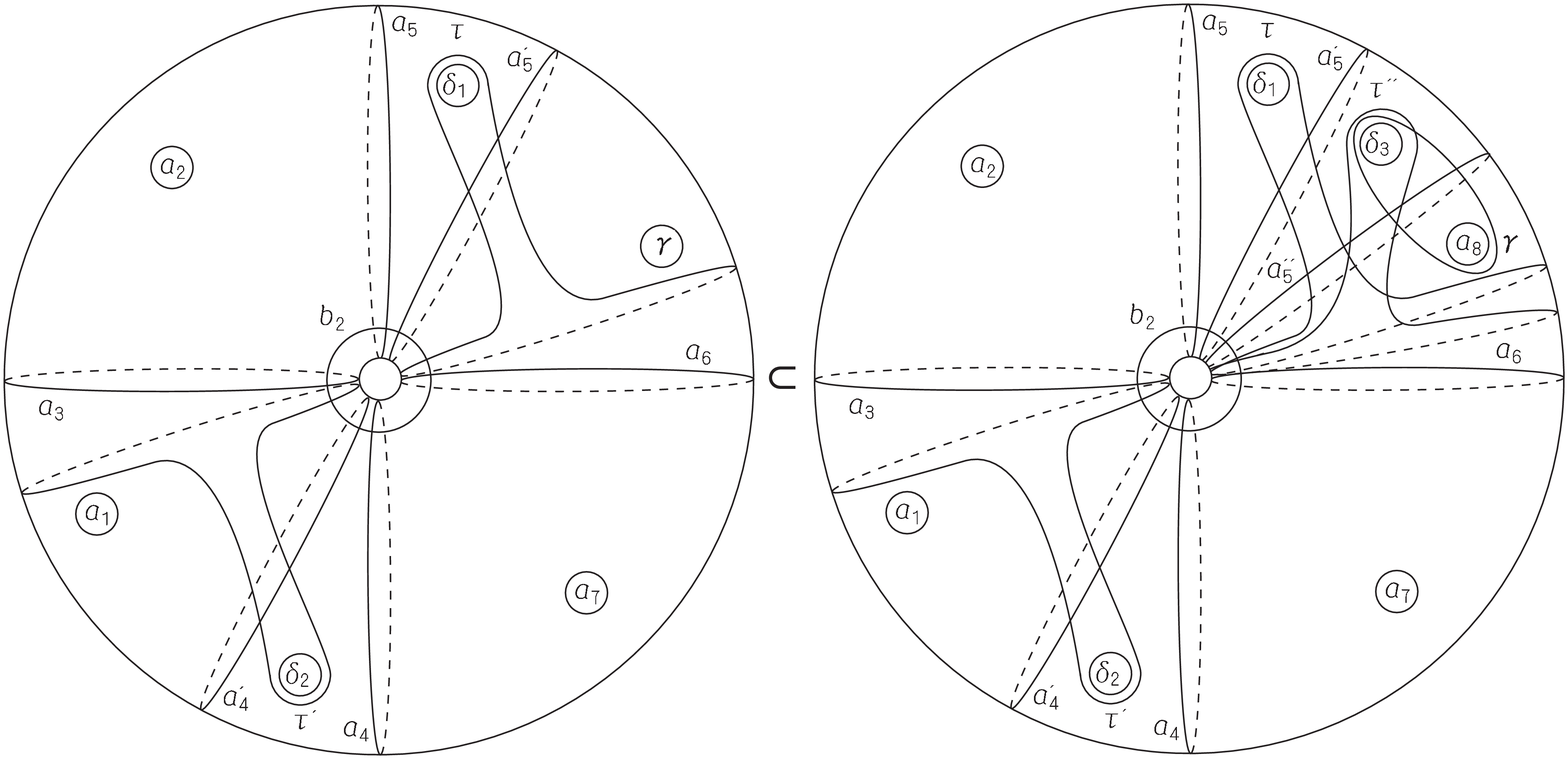}
\caption{\label{6to7'} Embedding of $\Sigma_{1,6}$ into $\Sigma_{1,7}$ (II)}
\end{figure}

We change the name $a_8$ of a curve in the relation (B2) 
into $\gamma$ (shown on the right in Figure \ref{6to7'}) to obtain 
\[
a_1a_2a_7\gamma\delta_1\delta_2 
\equiv b_2a_6a_3b_2\bar{a}_3\tau' a'_4a_5b_2a_6a_3b_2\bar{a}_6\tau a'_5a_4. 
\]
We embed $\Sigma_{1,6}$ into $\Sigma_{1,7}$ and take simple closed curves 
$a''_5,a_8,\delta_3,\tau''$ in the interior of $\Sigma_{1,7}$ 
shown in Figure \ref{6to7'}. Then we have a lantern relation 
\[
\delta_3a_8a_6a'_5 \equiv \gamma\tau'' a''_5. 
\]
Combining these relations 
and applying commutativity relations and conjugations, we obtain a relation 
\begin{align*}\tag{B3}
a_1a_2a_7a_8\delta_1\delta_2\delta_3 
& \equiv \bar{a}_6b_2a_6a_3b_2\bar{a}_3\tau' a'_4a_5b_2a_6a_3b_2\bar{a}_6\tau a'_5a_4 
 \bar{\gamma}\cdot \bar{a}'_5 \gamma \tau'' a''_5 \\
& \equiv \bar{a}_6b_2a_6a_3b_2\bar{a}_3\tau' a'_4a_5b_2a_6a_3b_2\bar{a}_6 
\tau\tau'' a_4a''_5 \\          
& \equiv b_2a_3a_6b_2\bar{a}_6\tau\tau'' a_4a''_5\bar{a}_6b_2a_6a_3b_2 
\bar{a}_3\tau' a'_4a_5.        
\end{align*}

\begin{figure}[h]
\includegraphics[width=10cm]{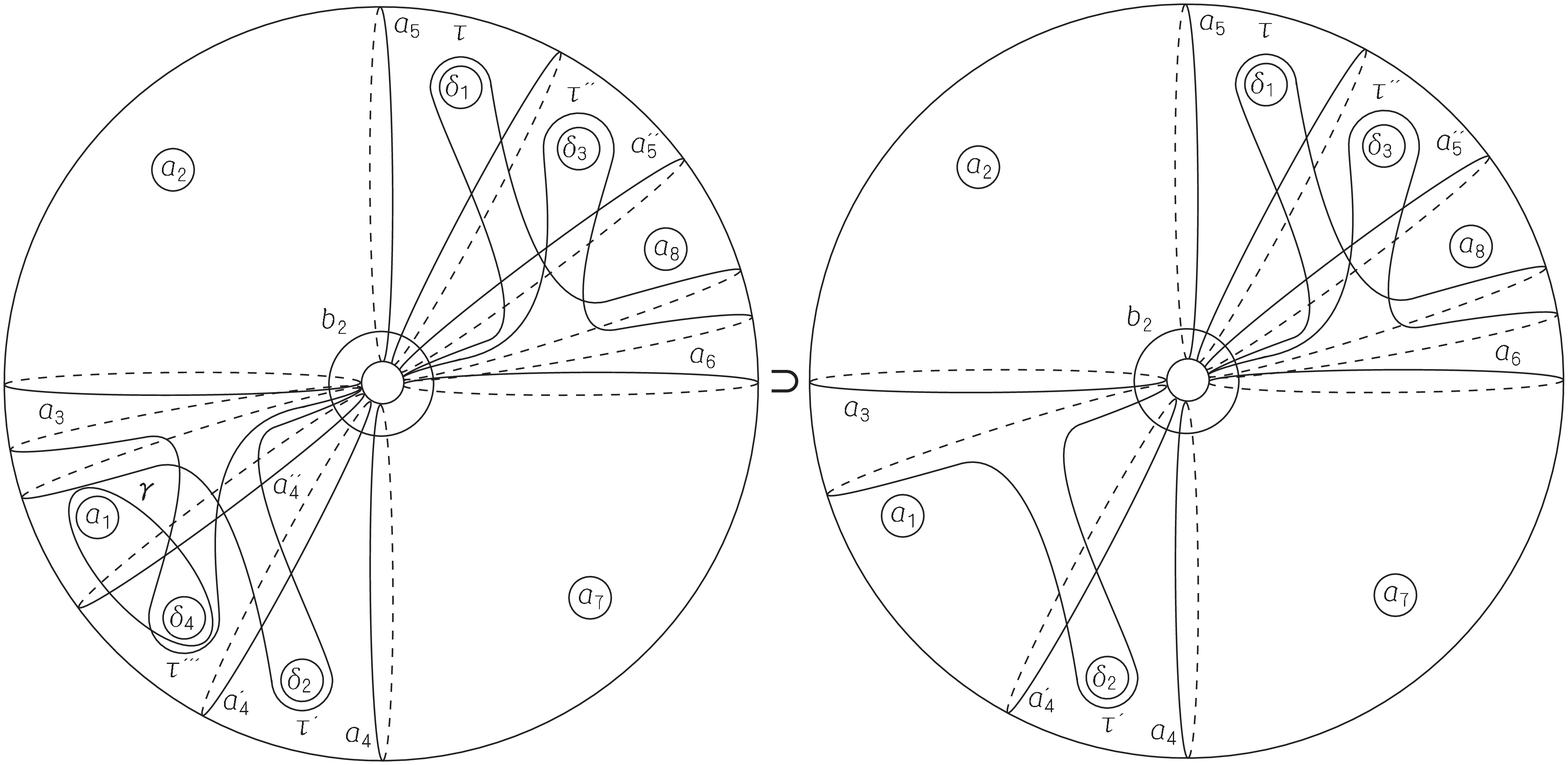}
\caption{\label{7to8'} Embedding of $\Sigma_{1,7}$ into $\Sigma_{1,8}$ (II)}
\end{figure}

We change the name $a_1$ of a curve in the relation (B3) 
into $\gamma$ (shown on the right in Figure \ref{7to8'}) to obtain 
\[
\gamma a_2a_7a_8\delta_1\delta_2\delta_3 
\equiv b_2a_3a_6b_2\bar{a}_6\tau\tau'' a_4a''_5 
\bar{a}_6b_2a_6a_3b_2\bar{a}_3\tau' a'_4a_5. 
\]
We embed $\Sigma_{1,7}$ into $\Sigma_{1,8}$ and take simple closed curves 
$a_1,a''_4,\delta_4,\tau'''$ in the interior of $\Sigma_{1,8}$ 
shown in Figure \ref{7to8'}. Then we have a lantern relation 
\[
\delta_4a_1a_3a'_4 \equiv \gamma\tau''' a''_4. 
\]
Combining these relations 
and applying commutativity relations and conjugations, 
we finally obtain Relation B: 
\begin{align*}                
a_1a_2a_7a_8\delta_1\delta_2\delta_3\delta_4 
& \equiv \bar{a}_3b_2a_3a_6b_2\bar{a}_6\tau\tau'' 
a_4a''_5\bar{a}_6b_2a_6a_3b_2\bar{a}_3\tau' a'_4a_5\bar{\gamma} 
\cdot \bar{a}'_4\gamma\tau''' a''_4 \\ 
& \equiv \bar{a}_3b_2a_3a_6b_2\bar{a}_6\tau\tau'' 
a_4a''_5\bar{a}_6b_2a_6a_3b_2\bar{a}_3\tau'\tau''' a_5a''_4 \\   
& \equiv a_4a''_5\bar{a}_6b_2a_6a_3b_2\bar{a}_3\tau'\tau''' 
a_5a''_4\bar{a}_3b_2a_3a_6b_2\bar{a}_6\tau\tau''. 
\end{align*}
This completes the proof of Proposition \ref{relationB}. $\square$

\begin{rem} 
Both of Relations A and B are different from the eight-holed torus relation of 
Korkmaz and Ozbagci \cite{KO} though the constructions are similar. 
\end{rem}

\section{Constructions}

In this section we construct a new relation on a compact oriented surface 
of genus $g$ with $4g+4$ boundary components 
by combining copies of Relations A and B obtained in the previous section.

\subsection{Higher genus} 

We assume $g\geq 3$. 
For integers $m,n\; (0<m\leq n)$ and 
words $W_m,W_{m+1},\ldots ,W_n\in\mathcal{F}_{g,r}$, we denote the product 
$W_mW_{m+1}\cdots W_n$ (resp. $W_n\cdots W_{m+1}W_m$) 
by $\prod_{i=m}^n W_i$ (resp. $\prod_{i=n}^m W_i$). 

\begin{figure}[h]
\includegraphics[width=11cm]{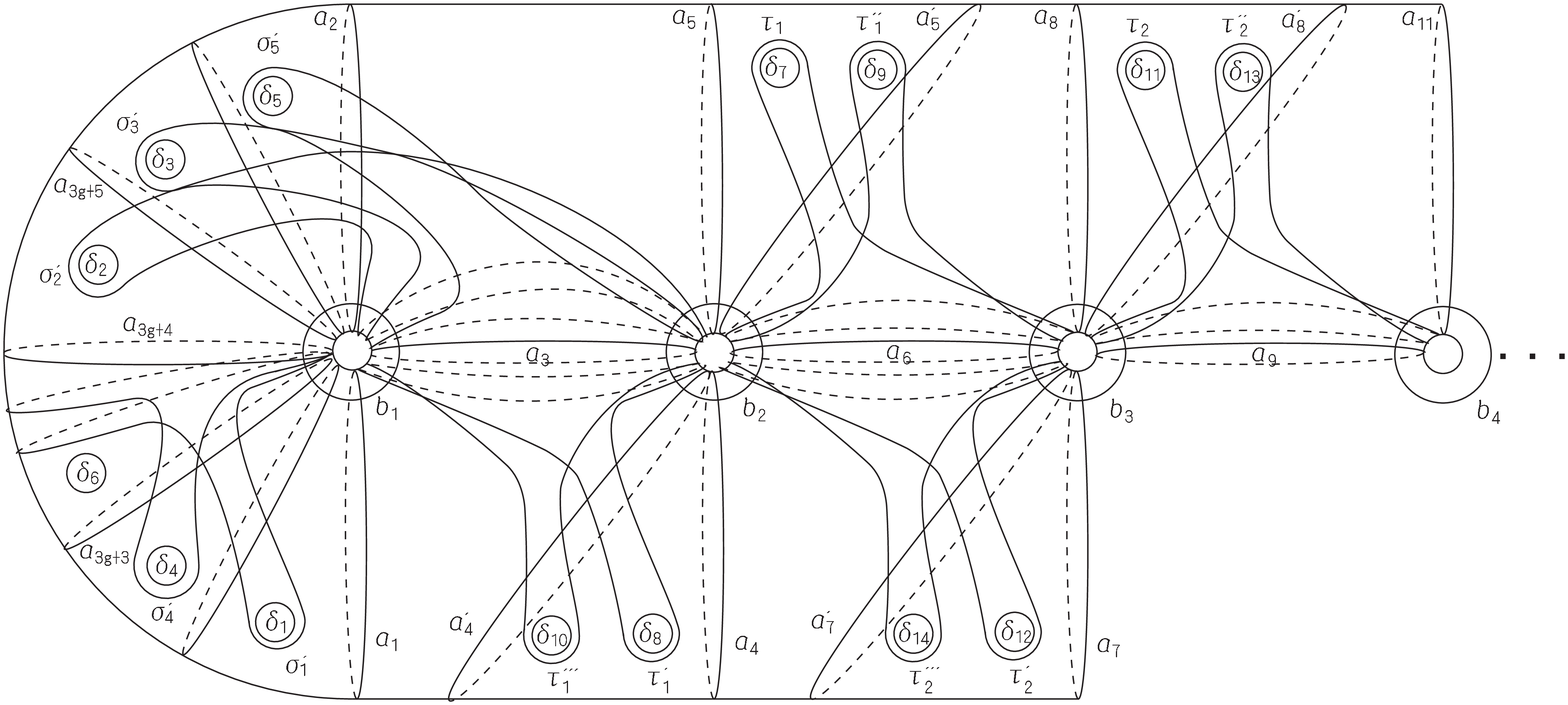}
\includegraphics[width=11cm]{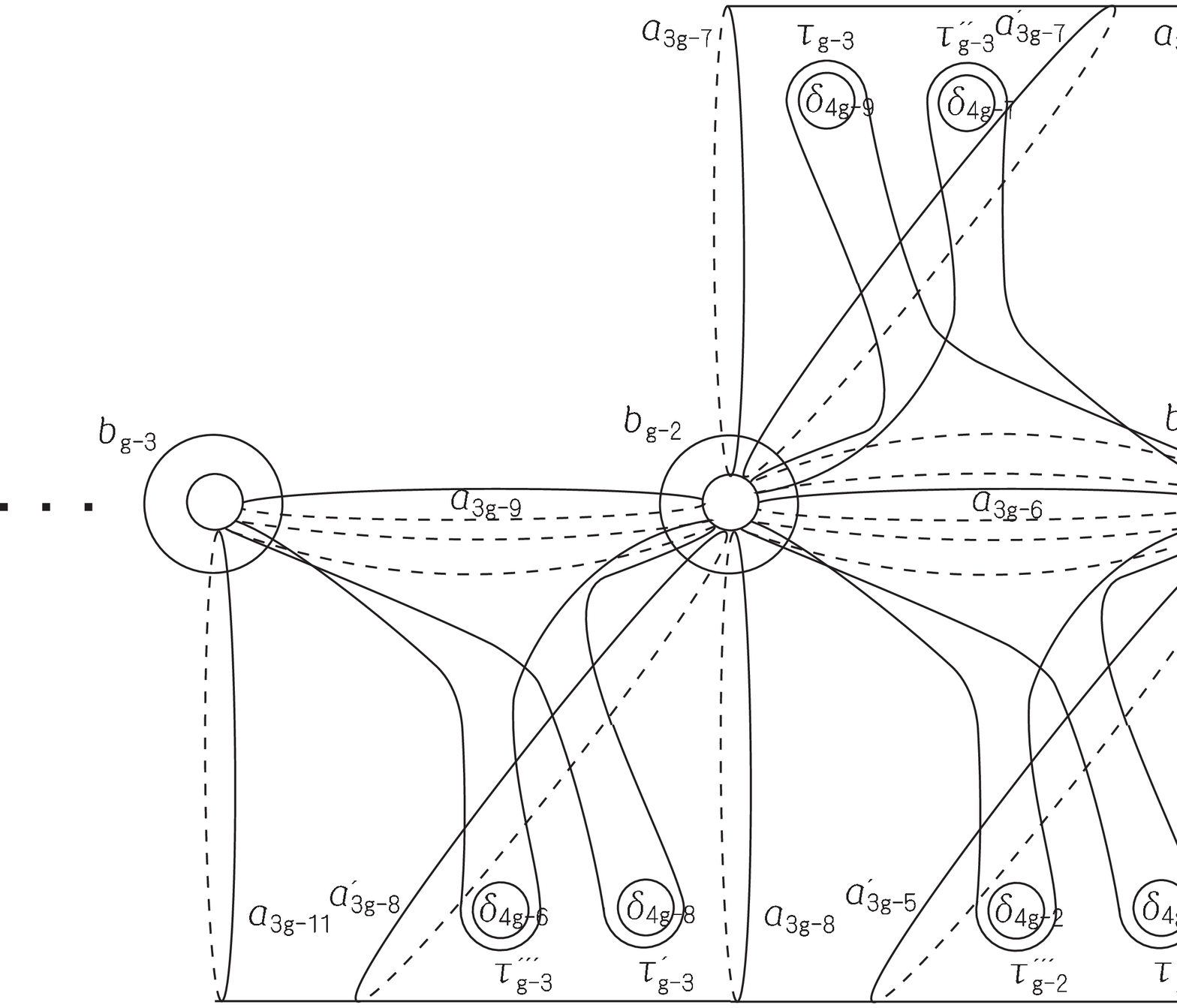}
\caption{\label{lift} Relation ${\rm H}_g$ for $g\geq 3$}
\end{figure}

\begin{thm}[Relation ${\rm H}_g$]\label{lift_hyp}
For simple closed curves in the interior of $\Sigma_{g,4g+4}$ shown in Figure \ref{lift}, 
we have the relation 
\begin{align*}
\delta_1\delta_2\dotsb\delta_{4g+3}\delta_{4g+4} 
& \equiv \prod_{i=g-1}^2 
\beta'''_i\beta_i\tau'_{i-1}\tau'''_{i-1} 
\cdot \beta_1\sigma'_1\sigma'_4a_{3g+3}\beta'_1\sigma'_2a_1\beta''_1\sigma'_3\sigma'_5 \\ 
& \quad\cdot \prod_{i=2}^{g-1} 
\beta''_i\beta'_i\tau_{i-1}\tau''_{i-1} 
\cdot \beta_g\sigma_1\sigma_4a_{3g}\beta'_g\sigma_2a_{3g-1}\beta''_g\sigma_3\sigma_5 
\end{align*}
in $\mathcal{M}_{g,4g+4}$, 
where $\beta_1:= {}_{a_{3g+4}}(b_1)$, $\beta'_1:={}_{a_3}(b_1)$, 
$\beta''_1:={}_{a_{3g+5}a_3}(b_1)$, $\beta_g:={}_{a_{3g+1}}(b_g)$, 
$\beta'_g:={}_{a_{3g-3}}(b_g)$, $\beta''_g:={}_{a_{3g-3}a_{3g+2}}(b_g)$, 
$\beta_i:={}_{a_{3i-3}}(b_i)$, $\beta'_i:={}_{a_{3i}}(b_i)$, 
$\beta''_i:={}_{\bar{a}_{3i-3}}(b_i)$ and $\beta'''_i:={}_{\bar{a}_{3i}}(b_i)$ 
$(i=2,\ldots ,g-1)$. 
\end{thm}

%\vspace{-0.7cm}

\begin{figure}[h]
\includegraphics[width=10cm]{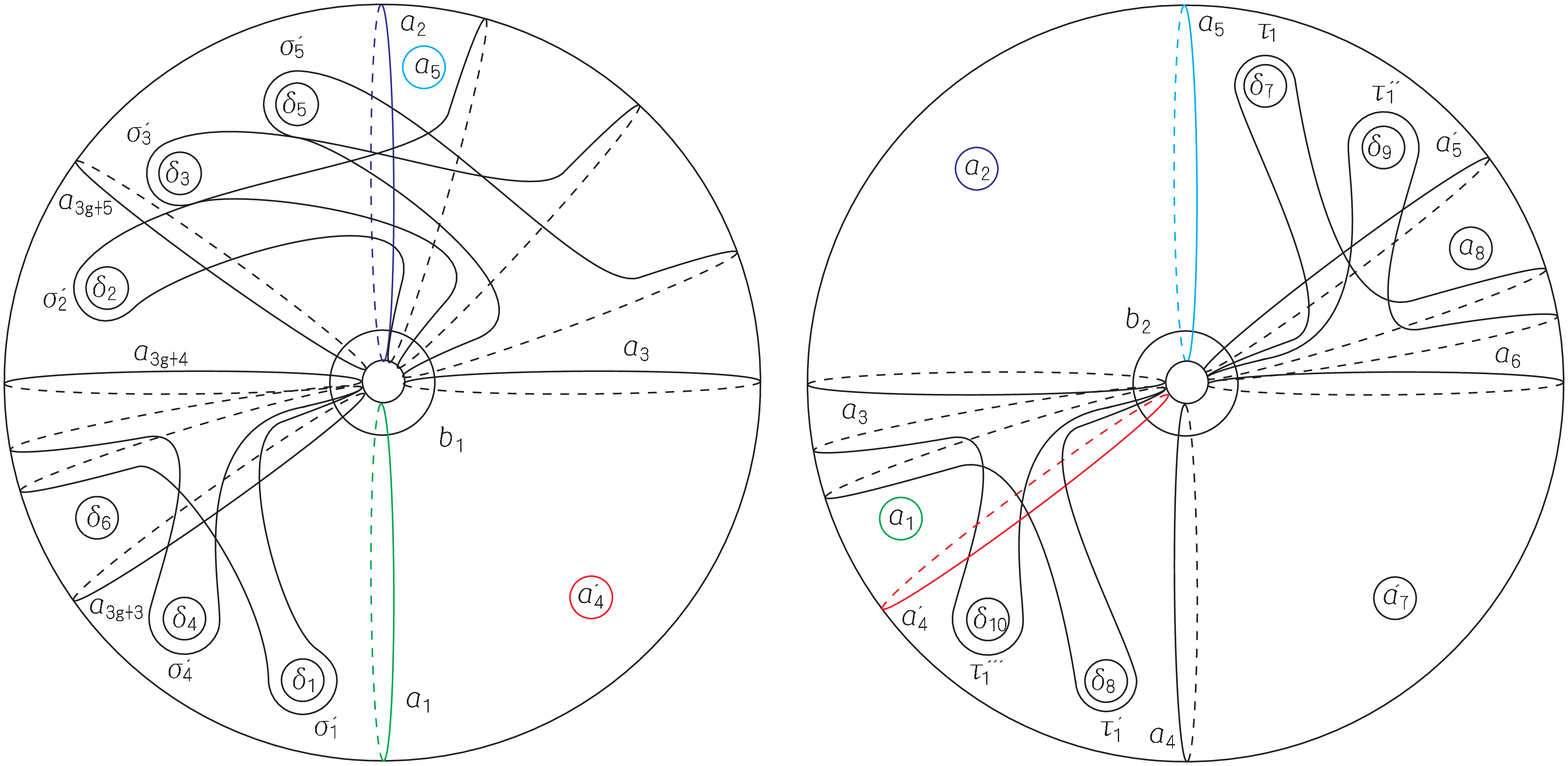}
\caption{\label{ab} Relations A and B}
\end{figure}

\noindent
{\it Proof}. We combine two copies of Relation A and $g-2$ copies of Relation B 
to obtain the desired relation. 
We first consider two relations for simple closed curves shown in Figure \ref{ab}. 
One is a copy of Relation A: 
{\allowdisplaybreaks %
\begin{align*}                
a'_4a_5\delta_1\delta_2\delta_3\delta_4\delta_5\delta_6 
& \equiv a_1a_2a_{3g+4}b_1\bar{a}_{3g+4}\sigma'_1\sigma'_4a_{3g+3}
a_3b_1\bar{a}_3\sigma'_2a_1a_3a_{3g+5}b_1\bar{a}_{3g+5}\bar{a}_3\sigma'_3\sigma'_5 \\
& \equiv a_1a_2\beta_1\sigma'_1\sigma'_4a_{3g+3}
\beta'_1\sigma'_2a_1\beta''_1\sigma'_3 \sigma'_5. 
\end{align*}}
Note that $\beta_1\equiv a_{3g+4}b_1\bar{a}_{3g+4}$, 
$\beta'_1\equiv a_3b_1\bar{a}_3$  
and $\beta''_1\equiv a_{3g+5}a_3b_1\bar{a}_3\bar{a}_{3g+5}$ by Remark \ref{braid}. 
The other is a copy of Relation B: 
\begin{align*}                
a_1a_2a'_7a_8\delta_7\delta_8\delta_9\delta_{10} 
& \equiv a_5a'_4 \bar{a}_3 b_2 a_3 a_6 b_2 \bar{a}_6 \tau_1 \tau''_1 
a_4 a'_5 \bar{a}_6 b_2 a_6 a_3 b_2 \bar{a}_3 \tau'_1 \tau'''_1. 
\end{align*}

\begin{figure}[h]
\includegraphics[width=8cm]{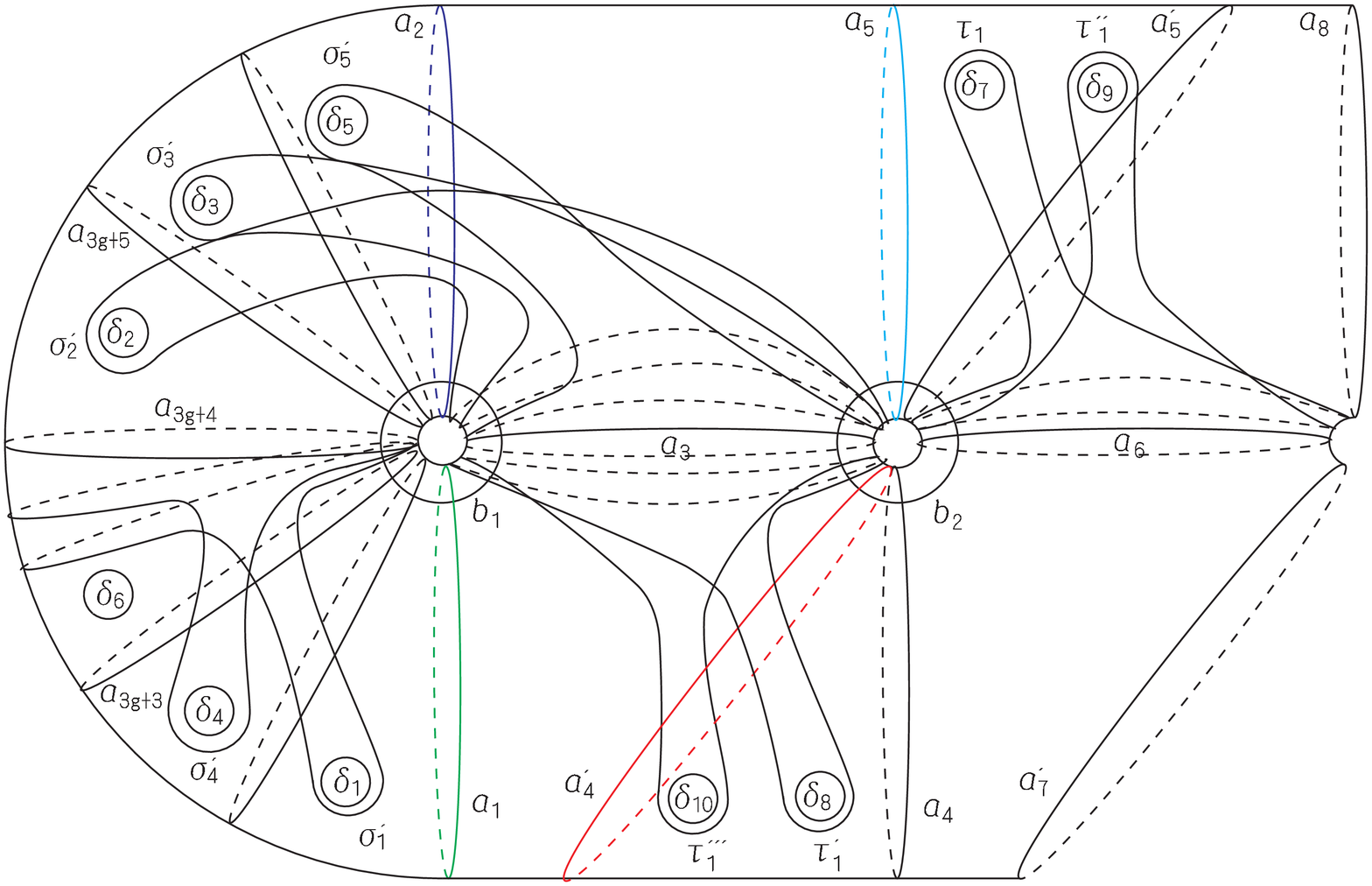}
\caption{\label{c2} Embeddings of two copies of 
$\Sigma_{1,8}$ into $\Sigma_{2,12}$ (I)}
\end{figure}

\noindent
We embed two copies of $\Sigma_{1,8}$ in Figure \ref{ab} 
into $\Sigma_{2,12}$ as shown in Figure \ref{c2}. 

%\vspace{-0.5cm}

\noindent
Combining these relations 
and applying commutativity relations and conjugations, 
we obtain a relation 
\begin{align*}\tag{C2}
& a'_7a_8\delta_1\delta_2\delta_3\delta_4\delta_5
\delta_6\delta_7\delta_8\delta_9\delta_{10} \\
\equiv & \; \bar{a}_5 \bar{a}'_4 a_5 a'_4 \bar{a}_3 b_2 a_3 a_6 b_2 \bar{a}_6 \tau_1 \tau''_1 
a_4 a'_5 \bar{a}_6 b_2 a_6 a_3 b_2 \bar{a}_3 \tau'_1 \tau'''_1 
\bar{a}_1 \bar{a}_2 \\
& \qquad \qquad \qquad \qquad \qquad \qquad \qquad 
\cdot a_1 a_2 \beta_1 \sigma'_1 \sigma'_4 a_{3g+3} 
\beta'_1 \sigma'_2 a_1 \beta''_1 \sigma'_3 \sigma'_5 \\
\equiv & \; \bar{a}_3 b_2 a_3 a_6 b_2 \bar{a}_6 \tau_1 \tau''_1a_4 a'_5 
\bar{a}_6 b_2 a_6 a_3 b_2 \bar{a}_3 \tau'_1 \tau'''_1   
\cdot \beta_1 \sigma'_1 \sigma'_4 a_{3g+3} \beta'_1 \sigma'_2 a_1 
\beta''_1 \sigma'_3 \sigma'_5 \\
\equiv & \; a_4 a'_5 \bar{a}_6 b_2 a_6 a_3 b_2 \bar{a}_3 \tau'_1 \tau'''_1 
\cdot \beta_1 \sigma'_1 \sigma'_4 a_{3g+3} \beta'_1 \sigma'_2 a_1 
\beta''_1 \sigma'_3 \sigma'_5
\cdot \bar{a}_3 b_2 a_3 a_6 b_2 \bar{a}_6 \tau_1 \tau''_1. 
\end{align*}

\begin{figure}[h]
\includegraphics[width=5cm]{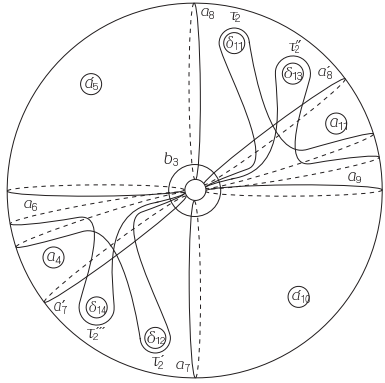}
\caption{\label{b3} Another Relation B}
\end{figure}

We next consider the relation (C2) and another copy of Relation B 
for simple closed curves shown in Figure \ref{b3}: 
\[
a_4 a'_5 a'_{10} a_{11} \delta_{11} \delta_{12} \delta_{13} \delta_{14} 
\equiv a_7 a'_8 \bar{a}_9 b_3 a_9 a_6 b_3 \bar{a}_6 \tau'_2 \tau'''_2 
a_8 a'_7 \bar{a}_6 b_3 a_6 a_9 b_3 \bar{a}_9 \tau_2 \tau''_2. 
\]

\begin{figure}[h]
\includegraphics[width=10cm]{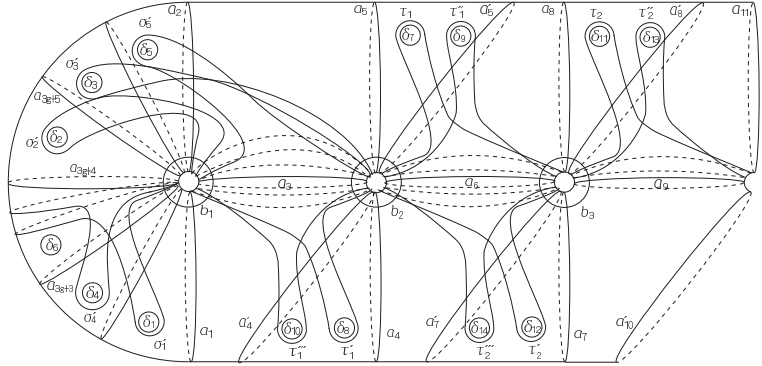}
\caption{\label{c3} Embeddings of $\Sigma_{2,12}$ and $\Sigma_{1,8}$ 
into $\Sigma_{3,16}$}
\end{figure}

\noindent
We embed $\Sigma_{2,12}$ in Figure \ref{c2} and $\Sigma_{1,8}$ in Figure \ref{b3} 
into $\Sigma_{3,16}$ as shown in Figure \ref{c3}. 

%\vspace{-0.5cm}

\noindent
Combining these relations 
and applying commutativity relations and conjugations, 
we obtain a relation 
{\allowdisplaybreaks %
\begin{align*}\tag{C3}
& a'_{10} a_{11} \delta_1 \delta_2 \delta_3 \delta_4 \delta_5 
\delta_6 \delta_{7} \delta_{8} \delta_{9} \delta_{10} \delta_{11} \delta_{12}  
\delta_{13} \delta_{14} \\
\equiv & \; \bar{a}'_7 \bar{a}_8 a_8 a'_7 \bar{a}_6 b_3 a_6 a_9 b_3 \bar{a}_9 \tau_2 \tau''_2 
a_7 a'_8 \bar{a}_9 b_3 a_9 a_6 b_3 \bar{a}_6 \tau'_2 \tau'''_2 \bar{a}_4 \bar{a}'_5 
\cdot a_4 a'_5 \bar{a}_6 b_2 a_6 a_3 b_2 \bar{a}_3 \tau'_1 \tau'''_1 \\
& \quad \cdot \beta_1 \sigma'_1 \sigma'_4 a_{3g+3} 
\beta'_1 \sigma'_2 a_1 \beta''_1 \sigma'_3 \sigma'_5
\cdot \bar{a}_3 b_2 a_3 a_6 b_2 \bar{a}_6 \tau_1 \tau''_1 \\ 
\equiv & \; \bar{a}_6 b_3 a_6 a_9 b_3 \bar{a}_9 \tau_2 \tau''_2 
a_7 a'_8 \bar{a}_9 b_3 a_9 a_6 b_3 \bar{a}_6 \tau'_2 \tau'''_2  
\cdot \bar{a}_6 b_2 a_6 a_3 b_2 \bar{a}_3 \tau'_1 \tau'''_1 \\
& \quad \cdot \beta_1 \sigma'_1 \sigma'_4 a_{3g+3} 
\beta'_1 \sigma'_2 a_1 \beta''_1 \sigma'_3 \sigma'_5
\cdot \bar{a}_3 b_2 a_3 a_6 b_2 \bar{a}_6 \tau_1 \tau''_1 \\ 
\equiv & \; a_7 a'_8 \bar{a}_9 b_3 a_9 a_6 b_3 \bar{a}_6 \tau'_2 \tau'''_2  
\cdot \bar{a}_6 b_2 a_6 a_3 b_2 \bar{a}_3 \tau'_1 \tau'''_1 \\
& \quad \cdot \beta_1 \sigma'_1 \sigma'_4 a_{3g+3} 
\beta'_1 \sigma'_2 a_1 \beta''_1 \sigma'_3 \sigma'_5
\cdot \bar{a}_3 b_2 a_3 a_6 b_2 \bar{a}_6 \tau_1 \tau''_1 
\cdot \bar{a}_6 b_3 a_6 a_9 b_3 \bar{a}_9 \tau_2 \tau''_2. 
\end{align*}}

\begin{figure}[h]
\includegraphics[width=11cm]{fig15.1.eps}
\includegraphics[width=11cm]{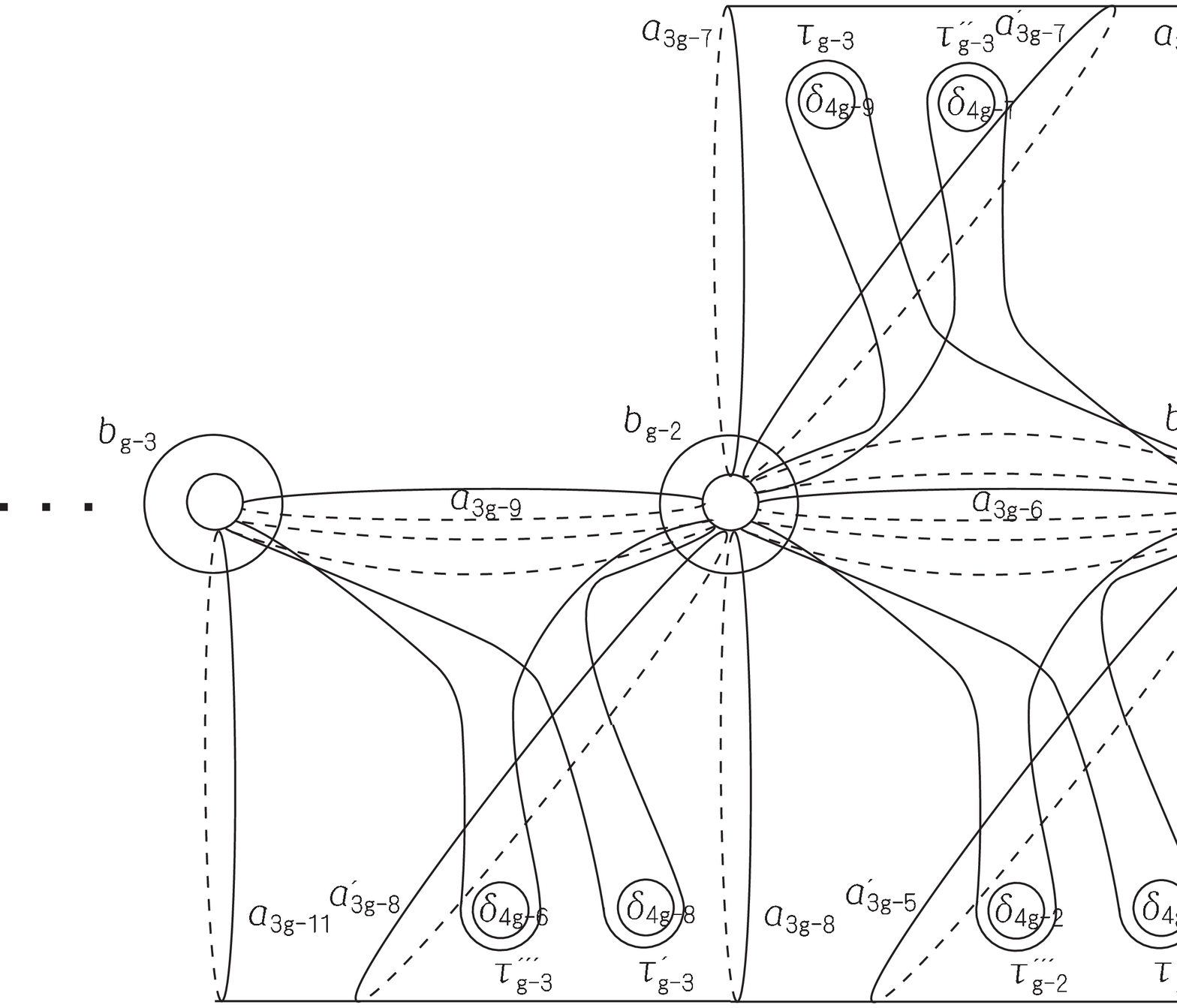}
\caption{\label{Cg-1} Relation (C($g-1$))}
\end{figure}

We repeat similar procedures by making use of $g-4$ copies of Relation B: 
{\allowdisplaybreaks %
\begin{align*}
& a_{3i-5}a'_{3i-4}a_{3i+1}a_{3i+2}\delta_{4i-1}\delta_{4i}\delta_{4i+1}\delta_{4i+2} \\
\equiv & \; a_{3i-2}a'_{3i-1}\bar{a}_{3i}b_ia_{3i}a_{3i-3}b_i\bar{a}_{3i-3}
\tau'_{i-1} \tau'''_{i-1} 
\cdot a_{3i-1}a'_{3i-2}\bar{a}_{3i-3}b_ia_{3i-3}a_{3i}b_i\bar{a}_{3i} 
\tau_{i-1} \tau''_{i-1}  
\end{align*}}
for $i=4,\ldots ,g-1$ to obtain relations (C4), (C5), $\cdots$, and 
\begin{align*}\tag{C($g-1$)}
& a_{3g-2} a_{3g-1}\delta_1\delta_2\cdots\delta_{4g-3}\delta_{4g-2} \\  
\equiv & \; a_{3g-5}a'_{3g-4} \prod_{i=g-1}^2 
\bar{a}_{3i}b_ia_{3i}a_{3i-3}b_i\bar{a}_{3i-3}\tau'_{i-1}\tau'''_{i-1} \\
& \quad 
\cdot \beta_1\sigma'_1\sigma'_4a_{3g+3}\beta'_1\sigma'_2a_1\beta''_1\sigma'_3\sigma'_5 
\cdot \prod_{i=2}^{g-1} 
\bar{a}_{3i-3}b_ia_{3i-3}a_{3i}b_i\bar{a}_{3i}\tau_{i-1}\tau''_{i-1} 
\end{align*}
for simple closed curves on $\Sigma_{g-1,4g}$ shown in Figure \ref{Cg-1}. 

%\vspace{-0.5cm}

\begin{figure}[h]
\includegraphics[width=5cm]{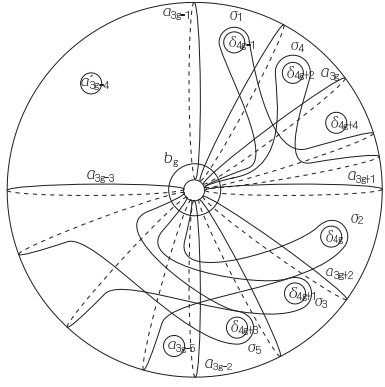}
\caption{\label{a} The other Relation A}
\end{figure}

We finally consider the other copy of Relation A for simple closed curves shown 
in Figure \ref{a}: 
{\allowdisplaybreaks %
\begin{align*}                
& a_{3g-5}a'_{3g-4}
\delta_{4g-1}\delta_{4g}\delta_{4g+1}\delta_{4g+2}\delta_{4g+3}\delta_{4g+4} \\
\equiv & \; a_{3g-2}a_{3g-1}\beta_g\sigma_1\sigma_4a_{3g} 
\beta'_g\sigma_2a_{3g-1}\beta''_g\sigma_3\sigma_5. 
\end{align*}}
Note that $\beta_g\equiv a_{3g+1}b_g\bar{a}_{3g+1}$, 
$\beta'_g\equiv a_{3g-3}b_g\bar{a}_{3g-3}$ and
$\beta''_g\equiv a_{3g-3}a_{3g+2}b_g\bar{a}_{3g+2}\bar{a}_{3g-3}$ 
by Remark \ref{braid}. 
We embed $\Sigma_{g-1,4g}$ in Figure \ref{Cg-1} and $\Sigma_{1,8}$ in Figure \ref{a} 
into $\Sigma_{g,4g+4}$ as shown in Figure \ref{lift}. 

\noindent
Combining these relations 
and applying commutativity relations and conjugations, 
we obtain Relation ${\rm H}_g$. 
Note that 
$\beta_i\equiv a_{3i-3}b_i\bar{a}_{3i-3}$, $\beta'_i\equiv a_{3i}b_i\bar{a}_{3i}$, 
$\beta''_i\equiv \bar{a}_{3i-3}b_ia_{3i-3}$ and $\beta'''_i\equiv \bar{a}_{3i}b_ia_{3i}$ 
$(i=2,\ldots ,g-1)$ by Remark \ref{braid}. 
Thus we complete the proof of Theorem \ref{lift_hyp}. 
$\square$

\subsection{Genus two} 

In this subsection we construct a relation on $\Sigma_{2,12}$ similar to relations 
constructed in the previous subsection. 

\begin{figure}[h]
\includegraphics[width=7.5cm]{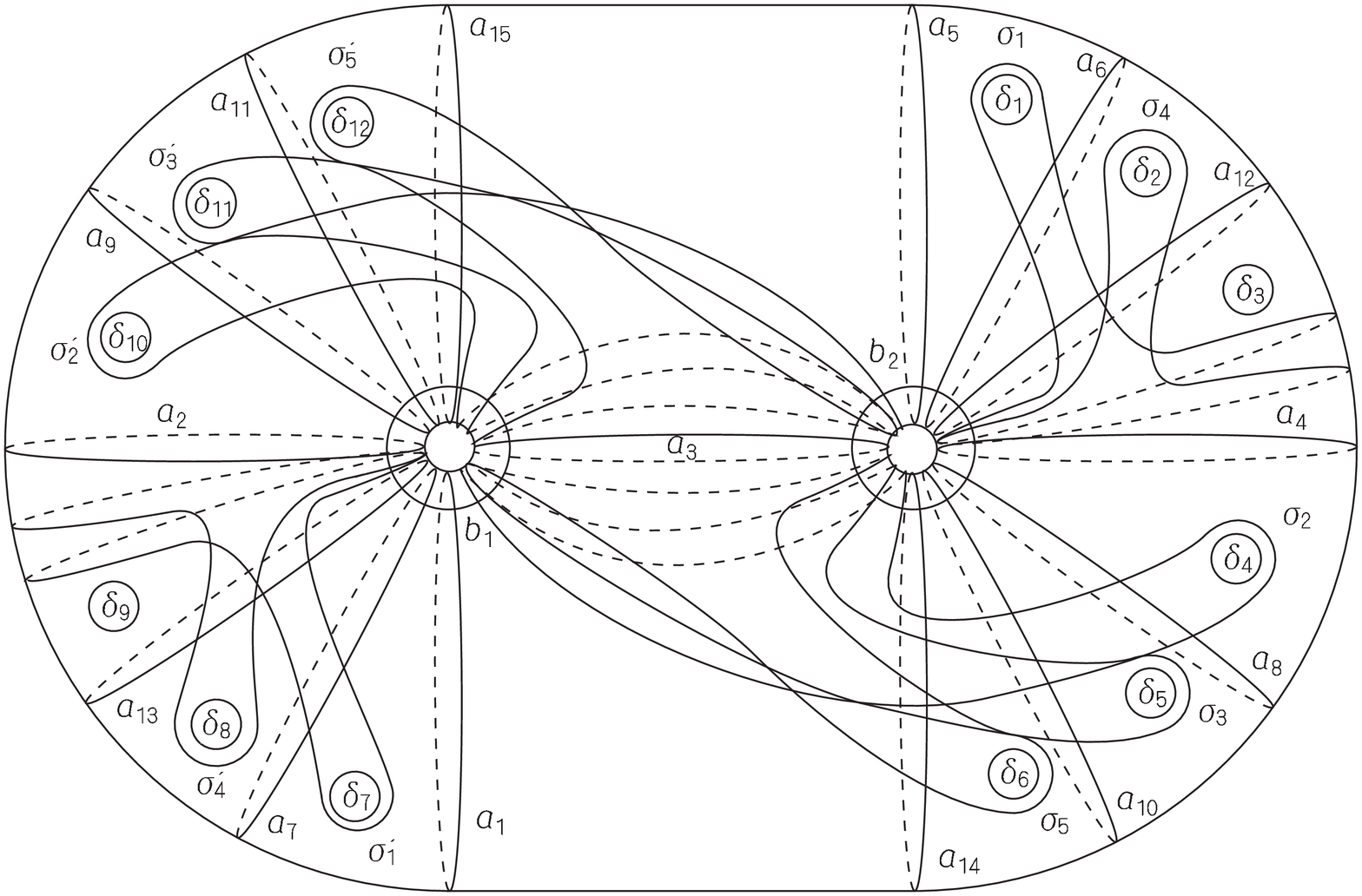}
\caption{\label{g=2} Embeddings of two copies of 
$\Sigma_{1,8}$ 
into $\Sigma_{2,12}$ (II)}
\end{figure}

\begin{thm}[Relation ${\rm H}_2$]\label{genus2}
For simple closed curves in the interior of $\Sigma_{2,12}$ shown in Figure \ref{g=2}, 
we have the relation 
\[
\delta_1 \delta_2 \delta_3 \delta_4 \delta_5 \delta_6 
\delta_7 \delta_8 \delta_9 \delta_{10} \delta_{11} \delta_{12} 
\equiv  \beta_1 \sigma'_1 \sigma'_4 a_{13} 
   \beta'_1 \sigma'_2 a_1 \beta''_1 \sigma'_3 \sigma'_5 
   \beta_2 \sigma_1 \sigma_4 a_{12} 
   \beta'_2 \sigma_2 a_5 \beta''_2 \sigma_3 \sigma_5 
\]
in $\mathcal{M}_{2,12}$, 
where $\beta_1:= {}_{a_2}(b_1)$, $\beta'_1:={}_{a_3}(b_1)$, 
$\beta''_1:={}_{a_3a_9}(b_1)$, $\beta_2:={}_{a_4}(b_2)$, 
$\beta'_2:={}_{a_3}(b_2)$ and $\beta''_2:={}_{a_3a_8}(b_2)$. 
\end{thm}

%\vspace{-0.5cm}

\begin{figure}[h]
\includegraphics[width=10cm]{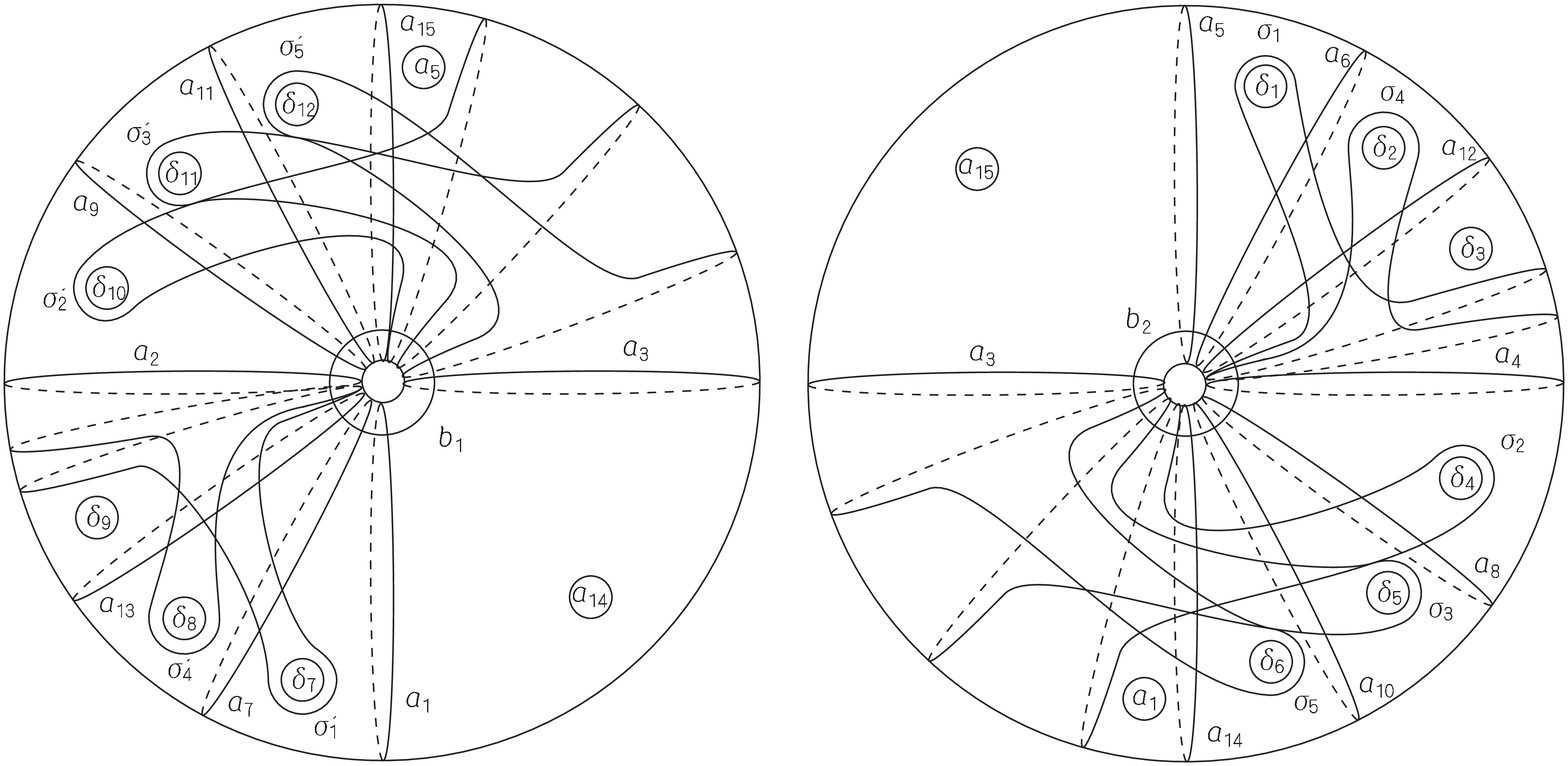}
\caption{\label{aa} Two copies of Relation A}
\end{figure}

\noindent
{\it Proof}. 
We first consider two copies of Relation A for simple closed curves shown in Figure \ref{aa}: 
\begin{align*}
a_5a_{14}\delta_7 \delta_8 \delta_9 \delta_{10}  \delta_{11} \delta_{12}
\equiv & \; a_{15} a_1 \beta_1 \sigma'_1 \sigma'_4 a_{13} 
\beta'_1 \sigma'_2 a_1 \beta''_1 \sigma'_3 \sigma'_5   \\
a_1 a_{15} \delta_1 \delta_2 \delta_3 \delta_4 \delta_5 \delta_6
\equiv & \; a_{14}  a_5 \beta_2 \sigma_1 \sigma_4 a_{12} 
\beta'_2 \sigma_2 a_5 \beta''_2 \sigma_3 \sigma_5. 
\end{align*}
Note that $ \beta_1 \equiv a_2 b_1 \bar{a}_2 $, $ \beta'_1 \equiv a_3 b_1 \bar{a}_3 $, 
$ \beta''_1 \equiv a_3 a_9 b_1 \bar{a}_9 \bar{a}_3 $,
$ \beta_2 \equiv a_4 b_2 \bar{a}_4 $, $ \beta'_2 \equiv a_3 b_2 \bar{a}_3 $  
and  $ \beta''_2 \equiv a_3 a_8 b_2 \bar{a}_8 \bar{a}_3 $ by Remark \ref{braid}. 

\noindent
Combining these relations 
and applying commutativity relations and conjugations, 
we obtain Relation ${\rm H}_2$. 
Thus we complete the proof of Theorem \ref{genus2}. 
$\square$

\section{Sections of Lefschetz fibrations}

In this section we show that the relation constructed in the previous section gives 
an explicit topological description of $4g+4$ disjoint $(-1)$-sections of 
a hyperelliptic Lefschetz fibration of genus $g$ on the manifold 
$\mathbb{CP}^2\# (4g+5)\overline{\mathbb{CP}}^2$. 

We begin with a definition of Lefschetz fibrations (cf. \cite{yukiomat}, \cite{GS}). 

\begin{defn}
Let $M$ be a closed oriented smooth $4$-manifold. A smooth map 
$f:M\rightarrow S^2$ is called a {\it Lefschetz fibration} of 
genus $g$ if it satisfies the following conditions: 

(i) $f$ has finitely many critical values $b_1,\ldots ,b_n\in 
S^2$ and $f$ is a smooth fiber bundle 
over $S^2-\{ b_1,\ldots ,b_n\}$ with fiber $\Sigma_{g,0}$; 

(ii) for each $i\; (i=1,\ldots ,n)$, there exists a unique 
critical point $p_i$ in the {\it singular fiber} $f^{-1}(b_i)$ 
such that $f$ is locally written as 
$f(z_1,z_2)=z_1^2+z_2^2$ with respect to some local complex 
coordinates around $p_i$ and $b_i$ which are compatible with 
orientations of $M$ and $S^2$; 

(iii) no fiber contains a $(-1)$-sphere.

\end{defn}

\begin{rem}\label{LF}
We always assume that a Lefschetz fibration is relatively minimal, it 
has at most one critical point on each fiber, and the genus of the base is 
equal to zero. A more general definition can be found in \cite[Chapter 8]{GS}. 
\end{rem}

Suppose that $g\geq 2$. 
According to theorems of Kas and Matsumoto, 
there exists a one-to-one correspondence between 
the isomorphism classes of Lefschetz fibrations and 
the equivalence classes of positive relators 
modulo simultaneous conjugations 
\[
c_1\cdot \cdots \cdot c_n \sim 
{}_W(c_1)\cdot \cdots \cdot {}_W(c_n), 
\]
and elementary transformations 
\begin{gather*}
c_1\cdot \cdots \cdot c_i\cdot c_{i+1}\cdot \cdots \cdot c_n \sim 
c_1\cdot \cdots \cdot c_{i+1}\cdot {}_{c_{i+1}^{-1}}(c_i)\cdot \cdots \cdot c_n, \\
c_1\cdot \cdots \cdot c_i\cdot c_{i+1}\cdot \cdots \cdot c_n \sim 
c_1\cdot \cdots \cdot {}_{c_i}(c_{i+1})\cdot c_i\cdot \cdots \cdot c_n, 
\end{gather*}
where $c_1\cdots c_n\in \mathcal{R}_{g,0}$ is 
a {\it positive relator} in the generators $\mathcal{S}_{g,0}$ 
and $W\in \mathcal{F}_{g,0}$. 
This correspondence is described by using the holonomy 
(or monodromy) homomorphism 
induced by the classifying map of $f$ restricted on $S^2-\{ b_1,\ldots ,b_n \}$ 
(cf. \cite{yukiomat}, \cite{GS}). 

\begin{defn}
Let $f:M\rightarrow S^2$ be a Lefschetz fibration of genus $g$. 
A smooth map $s:S^2\rightarrow M$ is called a {\it section} of $f$ 
if it satisfies $f\circ s={\rm id}_{S^2}$. A section $s$ of $f$ is an embedding of $S^2$ 
into $M$. The self-intersection number of the homology class 
$s_*([S^2])\in H_2(M;\mathbb{Z})$ is called the {\it self-intersection number} of $s$. 
A section of $f$ with self-intersection number $k$ is often called a $k$-{\it section}. 
\end{defn}

For a positive integer $r$, we attach $r$ disks 
to the boundary components of $\Sigma_{g,r}$ 
to obtain a closed surface $\Sigma_{g,0}$ and an embedding 
$\Sigma_{g,r}\hookrightarrow\Sigma_{g,0}$. 
This embedding induces a natural commutative diagram 
\begin{equation*}
\begin{CD}
1 @>>> \mathcal{R}_{g,r} @>>> \mathcal{F}_{g,r} @>{\varpi}>> \mathcal{M}_{g,r} @>>> 1 \\
@. @V{\lambda}VV @V{\lambda}VV @VVV @. \\
1 @>>> \mathcal{R}_{g,0} @>>> \mathcal{F}_{g,0} @>{\varpi}>> \mathcal{M}_{g,0} @>>> \, 1,
\end{CD}
\end{equation*}
where two horizontal sequences are exact. 
If two words $W_1$ and $W_2$ in $\mathcal{F}_{g,r}$ satisfy $W_1\equiv W_2$, 
then we have $\lambda(W_1)\equiv \lambda(W_2)$ in $\mathcal{F}_{g,0}$. In this case 
we call the relation $W_1\equiv W_2$ is a {\it lift} of the relation 
$\lambda(W_1)\equiv \lambda(W_2)$. 

\begin{lem}[cf. \cite{ABKP}, \cite{EKKOS}, \cite{smith}]\label{section}
Let $f:M\rightarrow S^2$ be a Lefschetz fibration of genus $g$ and 
$c_1\cdots c_n\in\mathcal{R}_{g,0}$ a positive relator corresponding to $f$. 
Suppose that there exists a relation $a_1\cdots a_n\equiv \delta_1^{k_1}\cdots \delta_r^{k_r}\, 
(a_1,\ldots ,a_n \in\mathcal{S}_{g,r},\, k_1,\ldots ,k_r>0)$ in $\mathcal{F}_{g,r}$ 
which is a lift of the relation $c_1\cdots c_n\equiv 1$ in $\mathcal{F}_{g,0}$, 
where $\delta_1,\ldots ,\delta_r$ are simple closed curves 
parallel to the boundary components of $\Sigma_{g,r}$. 
Then $f$ admits disjoint $r$ sections $s_1,\ldots ,s_r:S^2\rightarrow M$ with 
self-intersection number $-k_1,\ldots ,-k_r$, respectively. 
\end{lem}

For a chain $(c_1,\ldots ,c_{2g+1})$ of length $2g+1$ on $\Sigma_{g,0}$, 
we obtain a Lefschetz fibration $X_g\rightarrow S^2$ 
of genus $g$ associated to the hyperelliptic 
relation $(c_1\cdots c_{2g+1}c_{2g+1}\cdots c_1)^2$ $\equiv 1$ in 
$\mathcal{F}_{g,0}$. 
The total space $X_g$ of this fibration is known to be diffeomorphic to 
$\mathbb{CP}^2\# (4g+5)\overline{\mathbb{CP}}^2$ (cf. \cite{GS}, \cite{ito}). 

We denote the positive word on 
the right-hand side of Relation ${\rm H}_g$ by $U_g$ for $g\geq 2$. 
We consider the above embedding $\Sigma_{g,r}\hookrightarrow\Sigma_{g,0}$ 
and the commutative diagram for $r=4g+4$. 
By Theorems \ref{lift_hyp} and \ref{genus2}, 
Relation ${\rm H}_g$: $U_g\equiv \delta_1\delta_2\dotsb\delta_{4g+3}\delta_{4g+4}$ 
in $\mathcal{F}_{g,4g+4}$ is a lift of the relation $\lambda(U_g)\equiv 1$ in 
$\mathcal{F}_{g,0}$. 
This implies that the Lefschetz fibration $Y_g\rightarrow S^2$ of genus $g$ associated to 
the relation $\lambda(U_g)\equiv 1$ admits disjoint $4g+4$ sections with self-intersection 
number $-1$ by virtue of Lemma \ref{section}. 

\begin{thm}\label{isom}
Two Lefschetz fibrations $X_g$ and $Y_g$ are isomorphic to each other. 
\end{thm}

\noindent
{\it Proof}. Suppose that $g\geq 3$. 
We set $c_1:=\lambda(a_1),\,  c_{2g+1}:=\lambda(a_{3g-1}),\, 
c_{2i}:=\lambda(b_i)\, (i=1,\ldots ,g), \, c_{2i+1}:=\lambda(a_{3i})\, (i=1,\ldots ,g-1)$. 
Since $(a_1,b_1,a_3,b_2,\ldots ,a_{3g-3},b_g,a_{3g-1})$ is a chain of length 
$2g+1$ on $\Sigma_{g,4g+4}$, 
$(c_1,c_2,c_3,c_4,\ldots ,c_{2g-1},c_{2g},c_{2g+1})$ is a chain of length $2g+1$ 
on $\Sigma_{g,0}$. 
It is easily seen from Figure \ref{lift} that 
$\lambda(a_{3g+3})=\lambda(a_{3g+4})=\lambda(a_{3g+5})
=\lambda(\sigma'_1)=\lambda(\sigma'_4)=c_1$, 
$\lambda(a_{3g})=\lambda(a_{3g+1})=\lambda(a_{3g+2})
=\lambda(\sigma_1)=\lambda(\sigma_4)=c_{2g+1}$, 
$\lambda(\sigma'_2)=\lambda(\sigma'_3)=\lambda(\sigma'_5)=c_3$, 
$\lambda(\sigma_2)=\lambda(\sigma_3)=\lambda(\sigma_5)=c_{2g-1}$, 
$\lambda(\tau'_{i-1})=\lambda(\tau'''_{i-1})=c_{2i-1}$, 
$\lambda(\tau_{i-1})=\lambda(\tau''_{i-1})=c_{2i+1}$ $(i=2,\ldots ,g-1)$. 
Hence we obtain 
{\allowdisplaybreaks %
\begin{align*}
\lambda(U_g) & =
\prod_{i=g-1}^2 
({}_{\bar{c}_{2i+1}}(c_{2i}){}_{c_{2i-1}}(c_{2i})\cdot c_{2i-1}^2) \cdot 
{}_{c_1}(c_2)\cdot c_1^3\cdot {}_{c_3}(c_2)\cdot c_3c_1\cdot {}_{c_1c_3}(c_2)\cdot c_3^2 \\
& \quad\cdot \prod_{i=2}^{g-1} 
({}_{\bar{c}_{2i-1}}(c_{2i}){}_{c_{2i+1}}(c_{2i})\cdot c_{2i+1}^2) \\
& \quad\cdot {}_{c_{2g+1}}(c_{2g})\cdot c_{2g+1}^3\cdot 
{}_{c_{2g-1}}(c_{2g})\cdot c_{2g-1}c_{2g+1}
\cdot {}_{c_{2g-1}c_{2g+1}}(c_{2g})\cdot c_{2g-1}^2.
\end{align*}}
We now prove that $\lambda(U_g)\sim (c_1\cdots c_{2g+1}c_{2g+1}\cdots c_1)^2$ 
for $g\geq 3$. 
Applying elementary transformations (including cyclic permutations), 
we obtain the following sequence of equivalences. 
{\allowdisplaybreaks %
\begin{align*}
& \; \lambda(U_g) \\
\sim 
& \; c_{2g-1}\cdot \prod_{i=g-1}^2 
({}_{\bar{c}_{2i+1}}(c_{2i})\cdot c_{2i-1}c_{2i}c_{2i-1}) \cdot 
c_1c_2c_1^2c_3c_2c_1c_3\cdot {}_{c_1}(c_2)\cdot c_3 \\
& \cdot \prod_{i=2}^{g-1} 
({}_{\bar{c}_{2i-1}}(c_{2i})\cdot c_{2i+1}c_{2i}c_{2i+1}) 
\cdot c_{2g+1}c_{2g}c_{2g+1}^2c_{2g-1}c_{2g}c_{2g+1}c_{2g-1}
\cdot {}_{c_{2g+1}}(c_{2g}) \\
\sim 
& \; \prod_{i=g-1}^2 
(c_{2i+1}\cdot {}_{\bar{c}_{2i+1}}(c_{2i})\cdot c_{2i-1}c_{2i}) \cdot 
c_3c_1c_2c_1^2c_3c_2c_1c_3\cdot {}_{c_1}(c_2) \\
& \cdot \prod_{i=2}^{g-1} 
(c_{2i-1}\cdot {}_{\bar{c}_{2i-1}}(c_{2i})\cdot c_{2i+1}c_{2i}) 
\cdot c_{2g-1}c_{2g+1}c_{2g}c_{2g+1}^2c_{2g-1}c_{2g}c_{2g+1}c_{2g-1}
\cdot {}_{c_{2g+1}}(c_{2g}) \\
\sim 
& \; c_1c_{2g+1}\cdot \prod_{i=g-1}^2 
c_{2i}c_{2i+1}c_{2i-1}c_{2i} \cdot 
c_3c_2c_1^2c_3c_2c_1c_3\cdot {}_{c_1}(c_2) \\
& \cdot \prod_{i=2}^{g-1} 
c_{2i}c_{2i-1}c_{2i+1}c_{2i}
\cdot c_{2g-1}c_{2g}c_{2g+1}^2c_{2g-1}c_{2g}c_{2g+1}c_{2g-1}
\cdot {}_{c_{2g+1}}(c_{2g}) \\
\sim 
& \; \prod_{i=2g-2}^3 c_ic_{i+1} \cdot 
c_3c_2c_1^2c_3c_2c_1c_3\cdot {}_{c_1}(c_2)\cdot c_1 \\
& \cdot \prod_{i=3}^{2g-2} c_{i+1}c_i 
\cdot c_{2g-1}c_{2g}c_{2g+1}^2c_{2g-1}c_{2g}c_{2g+1}c_{2g-1}
\cdot {}_{c_{2g+1}}(c_{2g})\cdot c_{2g+1} \\
\sim 
& \; \prod_{i=2g-2}^3 c_ic_{i+1} \cdot 
c_3c_2c_3c_1^2c_2c_1^2c_3c_2 \\
& \cdot \prod_{i=3}^{2g-2} c_{i+1}c_i 
\cdot c_{2g-1}c_{2g}c_{2g-1}c_{2g+1}^2c_{2g}c_{2g+1}c_{2g-1}c_{2g+1}c_{2g} \\
\sim 
& \; \prod_{i=2g-2}^3 c_ic_{i+1} \cdot 
c_2c_3c_2c_1c_2c_1c_2c_1c_3c_2 \\
& \cdot \prod_{i=3}^{2g-2} c_{i+1}c_i 
\cdot c_{2g}c_{2g-1}c_{2g}c_{2g+1}c_{2g}c_{2g+1}c_{2g}c_{2g+1}c_{2g-1}c_{2g} \\
\sim 
& \; \prod_{i=2g-2}^2 c_ic_{i+1} \cdot 
c_1c_2c_1c_1c_2c_1 
\cdot \prod_{i=2}^{2g-1} c_{i+1}c_i 
\cdot c_{2g+1}c_{2g}c_{2g+1}c_{2g+1}c_{2g}c_{2g+1}c_{2g-1}c_{2g} \\
\sim 
& \; c_{2g}c_{2g+1}c_{2g-1}c_{2g}\cdot \prod_{i=2g-2}^1 c_ic_{i+1} \cdot 
c_1c_1 \cdot \prod_{i=1}^{2g} c_{i+1}c_i 
\cdot c_{2g+1}c_{2g+1} \\
\sim 
& \; \prod_{i=2g}^1 c_ic_{i+1} \cdot 
c_1c_1 \cdot \prod_{i=1}^{2g} c_{i+1}c_i 
\cdot c_{2g+1}c_{2g+1} 
\sim 
\; \prod_{i=2g+1}^1c_i\cdot \prod_{i=2g+1}^1c_i
\cdot \prod_{i=1}^{2g+1}c_i\cdot \prod_{i=1}^{2g+1}c_i \\
\sim 
& \; \prod_{i=1}^{2g+1}c_i \cdot \prod_{i=2g+1}^1c_i\cdot 
\prod_{i=1}^{2g+1}c_i\cdot \prod_{i=2g+1}^1c_i 
\; =\; (c_1\cdots c_{2g+1}c_{2g+1}\cdots c_1)^2.
\end{align*}}
Suppose that $g=2$. 
We set $c_1:=\lambda(a_1),\, c_2:=\lambda(b_1),\, 
c_3:=\lambda(a_3),\, c_4:=\lambda(b_2),\, c_5:=\lambda(a_5)$. 
Since $(a_1,b_1,a_3,b_2,a_5)$ is a chain of length 
$5$ on $\Sigma_{2,12}$, $(c_1,c_2,c_3,c_4,c_5)$ is a chain of length $5$ 
on $\Sigma_{2,0}$. It is easily seen from Figure \ref{g=2} that 
$\lambda(a_2)=\lambda(a_9)=\lambda(a_{13})
=\lambda(\sigma'_1)=\lambda(\sigma'_4)=c_1$, 
$\lambda(a_4)=\lambda(a_8)=\lambda(a_{12})
=\lambda(\sigma_1)=\lambda(\sigma_4)=c_5$, 
$\lambda(\sigma'_2)=\lambda(\sigma'_3)=\lambda(\sigma'_5)=
\lambda(\sigma_2)=\lambda(\sigma_3)=\lambda(\sigma_5)=c_3$. 
Hence we obtain 
{\allowdisplaybreaks %
\begin{align*}
\lambda(U_2)=
{}_{c_1}(c_2)\cdot c_1^3\cdot {}_{c_3}(c_2)\cdot c_3c_1\cdot {}_{c_3c_1}(c_2)
\cdot c_3^2\cdot {}_{c_5}(c_4)\cdot c_5^3\cdot {}_{c_3}(c_4)\cdot c_3c_5\cdot 
{}_{c_3c_5}(c_4)\cdot c_3^2. 
\end{align*}}
We now prove that $\lambda(U_2)\sim (c_1c_2c_3c_4c_5c_5c_4c_3c_2c_1)^2$. 
Applying elementary transformations (including cyclic permutations), 
we obtain the following sequence of equivalences. 
{\allowdisplaybreaks %
\begin{align*}
\lambda(U_2) 
& \sim \; c_1c_2c_1^2c_3c_2c_1c_3\cdot {}_{c_1}(c_2)\cdot c_3\cdot 
c_5c_4c_5^2c_3c_4c_5c_3\cdot {}_{c_5}(c_4)\cdot c_3 \\
& \sim \; c_5c_2c_1^2c_3c_2c_1c_3\cdot {}_{c_1}(c_2)\cdot c_3\cdot 
c_4c_5^2c_3c_4c_5c_3\cdot {}_{c_5}(c_4)\cdot c_3c_1 \\
& \sim \; c_2c_1^2c_3c_2c_1c_3\cdot {}_{c_1}(c_2)\cdot c_1c_3\cdot 
c_4c_5^2c_3c_4c_5c_3\cdot {}_{c_5}(c_4)\cdot c_3c_5 \\
& \sim \; c_2c_1^2c_3c_2c_1c_3c_1c_2c_3\cdot 
c_4c_5^2c_3c_4c_5c_3\cdot {}_{c_5}(c_4)\cdot c_5c_3 \\
& \sim \; c_2c_1^2c_3c_2c_1c_3c_1c_2c_3
c_4c_5^2c_3c_4c_5c_3c_5c_4c_3 \\
& \sim \; c_1c_3c_2c_3c_1c_1c_2c_3
c_4c_3c_5^2c_4c_5c_3c_5c_4c_3c_2c_1 \\
& \sim \; c_1c_2c_3c_2c_1c_1c_2c_4
c_3c_4c_5^2c_4c_5c_3c_5c_4c_3c_2c_1 \\
& \sim \; c_1c_2c_3c_4c_2c_1c_1c_2
c_3c_4c_5c_4c_5c_4c_3c_5c_4c_3c_2c_1 \\
& \sim \; c_1c_2c_3c_4c_2c_1c_1c_2
c_3c_5c_4c_5c_5c_4c_3c_5c_4c_3c_2c_1 \\
& \sim \; c_1c_2c_3c_4c_5c_2c_1c_1c_2
c_3c_4c_5c_5c_4c_3c_5c_4c_3c_2c_1 \\
& \sim \; c_1c_2c_3c_4c_5c_5c_4c_3c_5c_4c_3c_2c_1c_1c_2c_3c_4c_5c_2c_1 \\
& \sim \; c_1c_2c_3c_4c_5c_5c_4c_3c_2c_1c_1c_2c_3c_4c_5c_5c_4c_3c_2c_1 \\
& = \; (c_1c_2c_3c_4c_5c_5c_4c_3c_2c_1)^2. 
\end{align*}}
This completes the proof of Theorem \ref{isom}. 
$\square$

\medskip

The next corollary immediately follows from the theorem. 

\begin{cor}
The Lefschetz fibration $X_g\rightarrow S^2$ 
of genus $g$ associated to the hyperelliptic 
relation admits disjoint $4g+4$ sections with self-intersection 
number $-1$. 
\end{cor}

\begin{figure}[h]
\includegraphics[width=13cm]{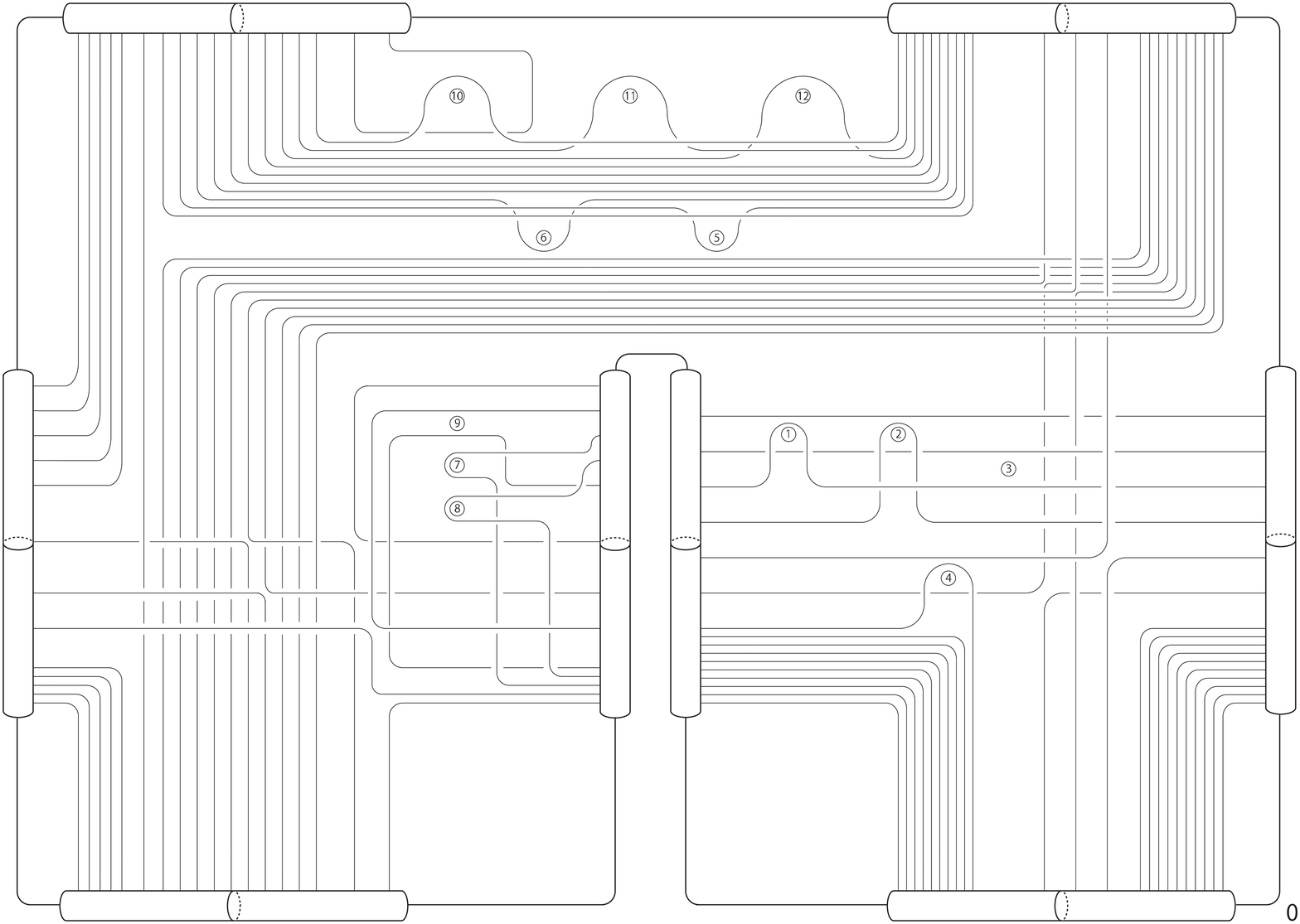}
\caption{\label{kirby} A Kirby diagram of $Y_2-\nu F$}
\end{figure}

By virtue of Theorem \ref{genus2}, we can even depict disjoint 
twelve sections of the Lefschetz fibration $Y_2\rightarrow S^2$ 
in a Kirby diagram of $Y_2-\nu F$, where $\nu F$ is an open 
fibered neighborhood of a regular fiber of $Y_2$ (cf. \cite[Sect. 4]{KO}). 
We first construct a handle decomposition of $\Sigma_{2,0}\times D^2$ with 
one $0$-handle, 
four $1$-handles, and one $2$-handle with framing $0$ from a fixed handle 
decomposition of $\Sigma_{2,0}$. 
We then attach twenty $2$-handles to $\Sigma_{2,0}\times D^2$ 
along the simple closed curves 
$\beta_1$, $\sigma'_1$, $\sigma'_4$, $a_{13}$, 
$\beta'_1$, $\sigma'_2$, $a_1$, $\beta''_1$, $\sigma'_3$, $\sigma'_5$, 
$\beta_2$, $\sigma_1$, $\sigma_4$, $a_{12}$, 
$\beta'_2$, $\sigma_2$, $a_5$, $\beta''_2$, $\sigma_3$, $\sigma_5$ (cf. Figure \ref{g=2}) 
on different fibers of $\Sigma_{2,0}\times S^1\rightarrow S^1$ 
with framing one less than the product framing of 
$\Sigma_{2,0}\times S^1$ 
to obtain a handle decomposition of $Y_2-\nu F$. 
Thus we have a Kirby diagram of $Y_2-\nu F$ shown in Figure \ref{kirby}. 
The framing coefficient of 
every component of the link but one with framing $0$ is equal to $-1$. 
Twelve disjoint sections coming from the simple closed curves 
$\delta_1,\ldots ,\delta_{12}$ are represented 
by twelve unknots transverse to each fiber of the fibration 
$\Sigma_{2,0}\times S^1\rightarrow S^1$ 
and meeting a fiber at twelve points indicated by encircled numbers 
$1,\ldots ,12$ in Figure \ref{kirby}. 
Attaching a $2$-handle with framing $-1$ along any one of the twelve unknots 
together with four $3$-handles and a $4$-handle to $Y_2-\nu F$, 
we have a handle decomposition of the closed manifold $Y_2$. 

By virtue of Theorem \ref{lift_hyp},
we can also depict disjoint 
$4g+4$ sections of the Lefschetz fibration $Y_g\rightarrow S^2$ 
in a Kirby diagram of $Y_g-\nu F$ for $g\geq 3$ 
in a similar way. 

The following proposition implies that the largest possible number of 
disjoint $(-1)$-sections of $X_g\rightarrow S^2$ is equal to $4g+4$ for most $g$. 

\begin{prop}\label{isom}
If $g$ is not equal to $k^2+k-1$ for any positive integer $k$, 
then the Lefschetz fibration $X_g\rightarrow S^2$ 
cannot admit disjoint $4g+5$ sections with self-intersection number $-1$. 
\end{prop}

\noindent
{\it Proof}. Suppose that the Lefschetz fibration $X_g\rightarrow S^2$ admits 
disjoint $4g+5$ sections $s_1,\ldots ,s_{4g+5}$ with self-intersection number $-1$. 
The orientation of $S^2$ induces that of $S_i:=s_i(S^2)$ for $i=1,\ldots ,4g+5$. 
We orient a regular fiber $F$ of $X_g$ so that it satisfies $[F]\cdot [S_i]=+1$ 
for $i=1,\ldots ,4g+5$. 
Blowing down the $(-1)$-spheres $S_1,\ldots ,S_{4g+5}$ in $X_g$, 
we obtain a $4$-manifold $X'$ and the image $F'$ of $F$ under the projection 
$X_g\rightarrow X'$. 
Since $[F]=[F']-[S_1]-\cdots -[S_{4g+5}]$ in $H_2(X_g;\mathbb{Z})\cong 
H_2(X';\mathbb{Z})\oplus (4g+5)H_2(\overline{\mathbb{CP}}^2;\mathbb{Z})$ and 
$[F]^2=0$, we have $[F']^2=4g+5$. 
On the other hand, $[F']^2$ must be the square of an integer because 
$[F']$ is a multiple of a generator of $H_2(X';\mathbb{Z})\cong \mathbb{Z}$. 
It is easy to see that $4g+5$ is the square of an integer if and only if 
$g$ is equal to $k^2+k-1$ for some positive integer $k$. 
$\square$

\begin{rem}\label{sato}
Two generic degree $d$ curves in $\mathbb{CP}^2$ induce a Lefschetz pencil of 
genus $(d-1)(d-2)/2$. Blowing up the base locus, we obtain a Lefschetz fibration 
$M_d\rightarrow S^2$ of the same genus. 
This fibration has $d^2$ sections with self-intersection number $-1$ and 
the total space $M_d$ is diffeomorphic to 
$\mathbb{CP}^2\# d^2\overline{\mathbb{CP}}^2$. It is well-known that 
the fibration $M_3\rightarrow S^2$ is isomorphic to $X_1\rightarrow S^2$, 
whereas the fibration $M_d\rightarrow S^2$ for $d\geq 4$ cannot be isomorphic to 
$X_g\rightarrow S^2$ for any $g$. 
\end{rem}

\medskip

\noindent
{\bf Acknowledgement.} 
The author is grateful to Professor Kazuhiro Konno for helpful comments 
on sections of fibrations on algebraic surfaces, and to Professor Yoshihisa Sato 
for helpful comments (Remark \ref{sato}) on Lefschetz fibrations.
And I wish to express my deep gratitude to my mentor, Hisaaki Endo,

\end{document}